\documentclass[12pt]{article}
\usepackage{latexsym, epsfig, amssymb, amsmath}
\usepackage[dvips]{color}
\usepackage{graphicx}
\usepackage{anysize}
\usepackage{hyperref}
\usepackage[comma,authoryear]{natbib}

\makeatletter

\renewcommand{\baselinestretch}{1.3}
\marginsize{1 in}{1 in}{1 in}{1 in}

\def\singlespace{\def\baselinestretch{1}\@normalsize}

\newcommand{\norm}[1]{\| #1 \|}

\newtheorem{theorem}{Theorem}

\newtheorem{lemma}{Lemma}

\def\Q{{\mathbf Q}}

\def\R{{\mathbf R}}
\def\P{{\mathbf P}}
\def\E{{\mathbf E}}
\def\1{{\mathbf 1}}        

\def\L{{\mathbf L}}
\def\I{{\mathbf I}}

\def\D{{\mathbf D}}
\def\H{\mathbf{H}}
\def\y{\mathbf{ y}}
\def\0{\mathbf{ 0}}

\def\W{\mathbf{W}}
\def\w{\mathbf{w}}
\def\A{\mathbf{A}}

\def\B{\mathbf{B}}

\def\x{\mathbf{x}}
\def\e{\mathbf{e}}

\def\M{\mathbf{M}}
\def\m{\mathbf{m}}

\def\a{\mathbf{a}}

\def\X{{\mathbf{X}}}

\def \I{\mathbf{I}}
\def \A{\mathbf{A}}

\def\bepsilon{\boldsymbol{\epsilon}}
\def\bSigma{\boldsymbol{\Sigma}}

\def\bl{\boldsymbol{\ell}}

\def\1{\mathbf{1}}
\def\0{\mathbf{0}}

\def\u{\mathbf{u}}

\def\tr{\text{tr}}

\def\blj0{\bl_j^{(0)}}

\def\diag{\text{diag}}
\def\var{\text{Var}}
\def\cov{\text{Cov}}
\def\sep{\text{sep}}

\makeatletter
\def\table{\@ifnextchar[{\table@i}{\table@i[\fps@table]}}
\def\table@i[#1]{\@float{table}[#1]\footnotesize}

\renewcommand{\hat}{\widehat}
\renewcommand{\tilde}{\widetilde}

\makeatother

\begin{document}

\title{\vspace{-0.9 in}
\bf  Estimation for Latent Factor Models for High-Dimensional Time
Series\thanks{Clifford Lam is Lecturer, Qiwei Yao is Professor, Department of Statistics, London
School of Economics, Houghton Street, London WC2A 2AE, U.K.
(email: C.Lam2@lse.ac.uk, Q.Yao@lse.ac.uk);
Neil Bathia is Research Fellow,
Department of Mathematics and Statistics, The University of
Melbourne, Victoria 3010 Australia (email:  nbathia@unimelb.edu.au).
Financial support from the  STICERD and LSE Annual Fund, and the Engineering and
Physical Sciences Research Council (UK) is
gratefully acknowledged.}
\date{}
\author{By Clifford Lam, Qiwei Yao and Neil Bathia  \\
Department of Statistics\\
London School of Economics and Political Science, and\\
Department of Mathematics and Statistics, University of Melbourne}} \maketitle

\vspace{-0.25 in}
\begin{singlespace}
\begin{quotation}
\indent
This paper deals with the dimension reduction for high-dimensional time series
based on common factors. In particular we allow the dimension of time series $p$ to be as
large as, or even larger than, the sample size $n$. The estimation for the factor loading
matrix and the factor process itself is carried out via an eigenanalysis for a $p\times p$
 non-negative definite matrix.
We show that when
all the factors are strong in the sense that the norm of each column
in the factor loading matrix is of the order $p^{1/2}$,
the estimator for the factor loading
matrix, as well as the resulting estimator for the precision matrix of the
original $p$-variant time series, are weakly consistent in $L_2$-norm with
the convergence rates independent of $p$. This result exhibits clearly
that the `curse' is canceled out by the `blessings' in
dimensionality.
We also establish the asymptotic properties of the estimation
when not all factors are strong. For the latter case, a two-step
estimation procedure is preferred accordingly to the asymptotic theory.
The proposed methods together with their asymptotic properties are further illustrated
in a simulation study. An application to a real data set is also reported.


\end{quotation}
\end{singlespace}

{\em Short Title}: Estimation of Large Factor Models.

{\em AMS 2000 subject classifications}. Primary 62F12, 62H25, 62H12.

{\em Key words and phrases}.
Convergence in $L_2$-norm,
curse and blessings of dimensionality,
dimension reduction,
eigenanalysis,
factor model,
precision matrix.

\newpage

\section{Introduction}\label{sect:introduction}
In this modern information age
analysis of large data sets is an integral part of both scientific research
and practical problem-solving.
High-dimensional time
series occur in many fields including, among others,
finance, economics, environmental and medical studies.
For example, to understand the dynamics of the returns of
large number of assets is the key for portfolio allocation, pricing and risk management.
Panel time series are common place in studying economic and business phenomena.
Environmental time series are often of a high-dimension because of the
large number of indices  monitored over many different locations.
On the other hand, the conventional time series models
such as vector AR or ARMA are not practically viable without a proper regularization when
the dimension is high, as the
number of parameters involved is a multiple of  the square of the dimension.
Hence it is pertinent to reduce the dimension of the
data before making further analysis. Different from the dimension-reduction
for independent observations, the challenge here is to retain the dynamical
structure of the time series.

Using common factors is one of the most frequently used and effective ways to
achieve dimension-reduction in analyzing multiple time series.
Early attempts in this direction include, for example,
\cite{A63}, \cite{PRT74}, \cite{B81} and \cite{PB87}.
To deal with the new challenge resulted from the fact that the number
of time series $p$ may be as large as, or even larger than, the length
of time series $n$ (such as most panel data),
more recent effort (mainly in econometrics) focuses on the inference when
$p$ goes to $\infty$ (along with
$n$). See, for example,  \cite{CR83}, \cite{C83}, \cite{B03},
\cite{FHLR00, FHLR04, FHLR02}. Furthermore motivated by analyzing some
economic and financial phenomena, those econometric factor models aim to identify
the common factors in the sense that each common factor affects the dynamics of
 {\sl most} of the original $p$ time series.  Those common factors are separated from
the so-called idiosyncratic `noise' components; each idiosyncratic component may
at most affect the dynamics of {\sl a few} original time series. Note an
idiosyncratic noise series is not necessarily white noise. The rigorous definition/identification
of the common factors and the idiosyncratic noise was established by \cite{CR83} and \cite{C83}
in an asymptotic manner when the number of time series goes to infinity, i.e. those
econometric factor models are only asymptotically identifiable when $p \to \infty$.
See also \cite{FHLR00}.


We adopt a different and more statistical approach in this paper from  a purely dimension-reduction
point of view. Our model is similar to those in \cite{PB87}, \cite{PP06} and \cite{PY08}.
However we consider the inference when $p$ is as large as, or even larger than, $n$.
Furthermore, we allow the future factors to depend on past (white) noise. This substantially
enlarge the capacity of the model.
Different from the aforementioned econometric factor models, we decompose the $p$-dimensional
time series into two parts: the dynamic part driven by a low-dimensional factor and the static
part which is a vector white noise.
Such a conceptually simple decomposition brings in conveniences in both model
identification and statistical inference.
In fact the model is identifiable for any finite $p$. Furthermore the
estimation for the factor loading matrix and the factor process itself is
equivalent to an eigenanalysis for a $p \times p$ non-negative definite matrix,
therefore is applicable when $p$ is in the order
of a few thousands. Our estimation procedure is rooted at the same idea as those
on which  the methods of \cite{PP06} and \cite{PY08} were based. However our method itself
is substantially simpler. For example, \cite{PP06} requires to compute the inverses of
sample autocovariance matrices,  which is computationally more costly when $p$ is large, and
is invalid when $p>n$. Furthermore in contrast to the eigenanalysis for one matrix,
it performs eigenanalysis for a matrix function of the sample autocovariance for several different
lags; see also \cite{PB87}.  The method of \cite{PY08} involves solving several nonlinear optimization
problems, which is designed to handle non-stationary factors and is only feasible for moderately large $p$.
Our approach identifies factors based on autocorrelation structure, which is more relevant
than the least squares approach advocated by \cite{BN02} and \cite{B03} in the context
of identifying time series factors.  In fact our method outperforms
the least squares method in a numerical experiment reported in section~\ref{sect:dataanalysis}.

The major theoretical contribution of this paper is to reveal an interesting and somehow intriguing
feature in factor modelling: the estimator for the factor loading matrix and the resulting estimator for
the precision matrix of the original $p$-dimensional time series converge
to the true ones at the rates independent of $p$,  provided that all the factors are strong in the
sense that the norm of each colunms in the factor loading matrix is of the order $p^{1/2}$.
Our simulation results indicate that indeed the estimation errors are indeed independent
of $p$. This result exhibits clearly that the `curse' is canceled out by the `blessings' in
dimensionality, as the high dimensionality is offset by combining
together the information from high-dimensional
data via common factors. However our factor model cannot improve
the estimation for the covariance matrix of the original time series, which coincides with the
result established by \cite{FFL08} with independent observations and known factors.

Another interesting finding from our asymptotic theory is to use a two-step estimation procedure
for a better performance when some factors are strong and some are not.
To this end, we characterize the strength of factors explicitly by an
index and show that the convergence rates of the estimators depend on
those indices. The concept of  weak and strong
factors was introduced in \cite{CPT09} in a different but related manner.

Further development of our setting with nonstationary factors together with forecasting issues
are reported in a companion paper \cite{LYB10b}.

The rest of the paper is organized as follows. The model, its presentational issues
and the estimation methods are presented in Section
\ref{sect:estimationmethod}.
Section \ref{sect:theories} contains the asymptotic properties of the
proposed methods when all the factors are of the same strength.
The results for the cases when there exist
factors of different levels of strength are given
in section \ref{sect:sswfactors}. Extensive simulation results are presented in section
\ref{sect:simulations}, with analysis of a set of implied volatility data
in section \ref{sect:dataanalysis}. All technical proofs are relegated to
section \ref{sect:Proofs}.

\section{Models and estimation methodology}\label{sect:estimationmethod}
\setcounter{equation}{0}
\subsection{Factor models}\label{subsect:factor_model}
{
Let $\y_1, \cdots, \y_n$ be $n$ $p\times 1$ successive observations from
a vector time series process.
}
The factor model assumes
\begin{equation}\label{eqn:factormodel}
  \y_t = \A\x_t + \bepsilon_t,
\end{equation}
where $\x_t$ is a $r \times 1$ unobserved factor time series which is assumed to be
strictly stationary with finite first two moments,
$\A$ is a $p \times r$ unknown constant
factor loadings matrix, and $r \leq p$ is the number of factors, and
$\{\bepsilon_t\}$ is a white noise with mean $\0$ and covariance matrix
$\bSigma_{\bepsilon}$. Furthermore, we assume that $\cov(\bepsilon_t, \x_s) = \0$ for all $s \leq t$,
and no linear combinations of the components of $\x_t$ are white noise. (Otherwise such combinations
should be absorbed in $\bepsilon_t$.)

Model (\ref{eqn:factormodel}) has been studied by, for example, \cite{PB87} and \cite{PP06} with
a stronger condition that the factor process and the white noise are uncorrelated
across all the lags. We relax this condition to allow the future factor
$\x_{t+k}$ correlated with the past white noise $\bepsilon_t$ ($k\ge 1$).
This is an appealing feature in modelling some economic and
finacial data.

%
%
%
%

In this paper, we always assume that the number of factors $r$ is known and fixed.
There is a large body of literature on how to determine $r$.
See, for example, \cite{BN02, BN07}, \cite{HL07}, \cite{PY08} and \cite{BYZ10}.
In section \ref{sect:simulations}, we use an information criterion proposed by \cite{BN02} to determine $r$
 in our simulation study.

\subsection{Identifiability and factor strength}\label{subsect:identifiability}
Model (\ref{eqn:factormodel}) is unchanged if
we replace the pair $(\A, \x_t)$ on the RHS by $(\A\H, \H^{-1}\x_t)$ for any invertible $\H$.
However the linear space spanned by the colunms of $\A$, denoted by
${\cal M}(\A)$ and called the factor loading space, is
uniquely defined by (\ref{eqn:factormodel}). Note ${\cal M}(\A)={\cal M}(\A \H)$
for any invertible $\H$.  Once such an $\A$ is specified, the factor process
$\x_t$ is uniquely defined accordingly.
We see the lack of uniqueness of $\A$ as an advantage, as we may choose a particular
$\A$ which facilitates our estimation in a simple and convenient manner.
Before we specify explicitly such an $\A$ in section~\ref{subsect:estimation}
below, we introduce an index $\delta$ for measuring the strength of factors, which is defined
naturally in terms of $\A$ in the first instance.
See conditions (A) and (B) below.


{
Let $\bSigma_{\x}(k) =\cov(\x_{t+k}, \x_t)$,\;\; $\bSigma_{\x,\bepsilon}(k)=
\cov(\x_{t+k},\; \bepsilon_t)$, and
\begin{align*}
 \tilde{\bSigma}_{\x}(k) &= (n-k)^{-1}\sum_{t=1}^{n-k} (\x_{t+k} - \bar{\x})(\x_t - \bar{\x})^T,\\
\tilde{\bSigma}_{\x,\bepsilon}(k) &= (n-k)^{-1}\sum_{t=1}^{n-k} (\x_{t+k} - \bar{\x})(\bepsilon_t - \bar{\bepsilon})^T,
 \end{align*}
with $\bar{\x} = n^{-1} \sum_{t=1}^n \x_t$, $\bar{\bepsilon} = n^{-1}\sum_{t=1}^n \bepsilon_t$. $\tilde{\bSigma}_{\bepsilon}(k)$ and $\tilde{\bSigma}_{\bepsilon,\x}(k)$ are defined similarly. Denote
by $\norm{\M}$ the spectral norm of $\M$, which is the positive square
root of the maximum eigenvalue of $\M\M^T$; and by $\norm{\M}_{\min}$
the positive square root of the minimum eigenvalue of $\M\M^T$ or
$\M^T\M$, whichever has a smaller matrix size. The notation $a \asymp
b$ represents $a = O(b)$ and $b = O(a)$.
}
Now we introduce the conditions
on the strength of factors.

\begin{itemize}
\item[(A)]
For $k=0,1,\cdots,k_0$, where $k_0 \geq 1$ is a small positive
integer, $\bSigma_{\x}(k)$ is full-ranked.
The cross autocovariance matrix $\bSigma_{\x,\bepsilon}(k)$ has elements of order $O(1)$.

\item[(B)] $\A = (\a_1 \cdots \a_r)$ such that
$\norm{\a_i}^2 \asymp p^{1-\delta}, \;\; i=1,\cdots,r, \;\; 0 \leq
\delta \leq 1.$

\item[(C)] For each
$i=1,\cdots,r$ { and $\delta$ given in (B)}, $\;\; \min_{\theta_j, j\neq i} \norm{ \a_i -
\sum_{j\neq i} \theta_j\a_j }^2 \asymp p^{1-\delta}$.
\end{itemize}

Note that model (\ref{eqn:factormodel}) is practically useful only if $r<<p$.
In our asymptotic theory we assume that $r$ remains as a constant while both $p$ and $n$
go to infinity. Therefore $\bSigma_{\x}(k)$ is an $r \times r$ fixed matrix of the full rank;
see (A).
%
When $\delta = 0$ in assumption (B), the corresponding factors are
called strong factors since it includes the case where each element
of $\a_i$ is $O(1)$, implying that the factors are shared (strongly) by
the majority of the $p$ cross-sectional variables. On the other
hand, they are called weak factors when $\delta > 0$.
\cite{CPT09} introduced a notion of strong
and weak factors, determined by  the finiteness of the mean absolute
values of the component of $\a_i$.
In this paper, we introduce index  $\delta$
which links explicitly the strength of factors $\x_t$ and the convergence rates
of our estimators. In fact the convergence is slower in the presence of
weak factors.
Assumptions (B) and (C) together ensure that all $r$ factors in the model
are of the equal strength.

To facilitate our estimation, we normalize the factor loadings matrix such that all
the columns of  $\A$ are orthonormal, i.e. $\A^T\A = \I_r$; see, e.g.
\cite{PY08}.
Then under assumptions (A) -- (C),  model (\ref{eqn:factormodel}) admits the follow
representation. Its proof is given  in the beginning of section~\ref{sect:Proofs} below.
\begin{equation}\label{eqn:factormodelWLOG}
\begin{split}
\y_t &= \A\x_t + \bepsilon_t, \;\; \A^T\A = \I_r, \;\text{ with }\\
\norm{\bSigma_{\x}(k)} &\asymp p^{1-\delta}
\asymp \norm{\bSigma_{\x}(k)}_{\min}, \;\;\; \norm{\bSigma_{\x,\bepsilon}(k)} = O(p^{1-\delta/2}),\\
\cov(\x_s, \bepsilon_t) &= \0 \;\; \text{ for all } s \leq t,\\
\cov(\bepsilon_s, \bepsilon_t) &= \0 \;\; \text{ for all } s \neq t.
\end{split}
\end{equation}

Unless specified otherwise, all $\y_t,\; \x_t$ and $ \bepsilon_t$
in the rest of this section and also section~\ref{sect:theories} are
defined in (\ref{eqn:factormodelWLOG}).

\subsection{Estimation}\label{subsect:estimation}

For $k \geq 1$, model (\ref{eqn:factormodelWLOG}) implies that
\begin{equation}
\bSigma_{\y}(k) = \cov(\y_{t+k}, \y_{t}) = \A\bSigma_{\x}(k)\A^T + \A\bSigma_{\x,\bepsilon}(k).
\end{equation}
For $k_0\ge 1$ given in condition (A),  define
\begin{equation} \label{eqn:Lexpansion}
\L = \sum_{k=1}^{k_0} \bSigma_{\y}(k)\bSigma_{\y}(k)^T
= \A \bigg( \sum_{k=1}^{k_0} \{\bSigma_{\x}(k)\A^T +
\bSigma_{\x,\bepsilon}(k)\}\{\bSigma_{\x}(k)\A^T +
\bSigma_{\x,\bepsilon}(k) \}^T\bigg)\A^T.
\end{equation}
Obviously $L$ is a $p\times p$ non-negative definite matrix.
Now we are ready to specify the factor loading matrix $\A$ to be used in our estimation.
First note that (\ref{eqn:factormodelWLOG}) is unchanged
if we replace $(\A, \x_t)$ by
$(\A \Q, \Q^T\x_t)$ for any $r\times r$ orthogonal matrix $\Q$.
Apply the spectrum decomposition to the positive-definite matrix
sandwiched by $\A$ and $\A^T$ on the RHS of (\ref{eqn:Lexpansion}), i.e.
$$ \sum_{k=1}^{k_0} \{\bSigma_{\x}(k)\A^T +
\bSigma_{\x,\bepsilon}(k)\}\{\bSigma_{\x}(k)\A^T +
\bSigma_{\x,\bepsilon}(k) \}^T = \Q \D \Q^T,
$$
where $\Q$ is an $r\times r$ orthogonal matrix, and $\D$ is a diagonal matrix with the elements
on the main diagonal in descending order.
This leads to
$
\L = \A \Q \D \Q^T \A^T.
$
As $\Q^T \A^T \A \Q = \I_r$, the columns of $\A \Q$ are the eigenvectors of $\L$ corresponding to
its $r$ non-zero eigenvalues. We take $\A \Q$ as the $\A $ to be used in our inference, i.e.
\begin{quote}
{\sl
the columns of the factor
loading matrix $\A $ are the  $r$  orthonormal eigenvectors of the matrix $\L$ corresponding
to its $r$ non-zero eigenvalues.
}
\end{quote}

A natural estimator for the $\A$ specified above is defined as
$\hat{\A} = (\hat{\a}_1, \cdots, \hat{\a}_r)$,
where $ \hat{\a}_i$ are the eigenvector of $\tilde{\L}$
corresponding to its $i$-th largest
eigenvalues, $\hat{\a}_1, \cdots, \hat{\a}_r$ are orthonormal, and
%
\begin{equation}\label{eqn:Ltilde}
  \tilde{\L} = \sum_{k=1}^{k_0} \tilde{\bSigma}_{\y}(k)
\tilde{\bSigma}_{\y}(k)^T, \;\;\; \tilde{\bSigma}_{\y}(k) =
{1 \over n- k}
\sum_{t=1}^{n-k} (\y_{t+k} - \bar{\y})(\y_t - \bar{\y})^T,
\end{equation}
with $\bar{\y} = n^{-1}\sum_{t=1}^n \y_t$.

Consequently, we estimate the factors and the residuals respectively by
\begin{equation}\label{eqn:estfactorresiduals}
  \hat{\x}_t = \hat{\A}^T\y_t, \;\;\;\quad\quad \e_t = \y_t - \hat{\A}\hat{\x}_t = (\I_p - \hat{\A}\hat{\A}^T)\y_t.
\end{equation}


\section{Asymptotic theory}
\label{sect:theories}
In this section we present the rates of convergence for the estimator
$\hat{\A}$ for model (\ref{eqn:factormodelWLOG}), as well as the
corresponding estimators for the covariance matrix and the precision
matrix, derived from model (\ref{eqn:factormodelWLOG}).
We need the following assumption for the original model (\ref{eqn:factormodel}):
\begin{itemize}
\item[(D)] { It holds for any $0\le k\le k_0$ } that
the elementwise rates of convergence for $\tilde{\bSigma}_{\x}(k) - \bSigma_{\x}(k)$, $\tilde{\bSigma}_{\x, \bepsilon}(k) - \bSigma_{\x,\bepsilon}(k)$ and $\tilde{\bSigma}_{\bepsilon}(k) - \bSigma_{\bepsilon}(k)$ are respectively $O_P(n^{-l_{x}})$, $O_P(n^{-l_{x\epsilon}})$ and $O_P(n^{-l_{\epsilon}})$, for some constants $0 < l_x, l_{x\epsilon}, l_{\epsilon} \leq 1/2$. We also have, elementwise, $\tilde{\bSigma}_{\bepsilon,\x}(k) = O_P(n^{-l_{x\epsilon}})$.
\end{itemize}
With the above assumption on the elementwise convergence for the sample cross- and auto-covariance matrices
of $\x_t$ and $\bepsilon_t$, we specify the convergence rate in
the spectral norm for the estimated factor loading matrix $\hat\A$.
It goes without saying explicitly that we may replace
some $\hat \a_j$ by $- \hat\a_j$ in order to match the direction of $\a_j$.

\begin{theorem}\label{thm:A}
Let assumption (D) hold, $\norm{\bSigma_{\x,\bepsilon}(k)} = o(p^{1-\delta})$, and
the $r$ non-zero eigenvalues of matrix $\L$ in (\ref{eqn:Lexpansion}) are different.
Then under model (\ref{eqn:factormodelWLOG}),
it holds that
$$ \norm{\hat{\A} - \A} = O_P(h_n) = O_P(n^{-l_{x}} +
p^{\delta/2}n^{-l_{x\epsilon}} + p^{\delta}n^{-l_{\epsilon}}), $$
provided
$h_n = o(1)$.
\end{theorem}
This theorem shows explicitly how the strength of the factors $\delta$
affects the rate of convergence. The convergence is faster
when the factors are stronger (i.e. $\delta$ gets smaller).
When $\delta = 0$, the rate is independent of $p$. This shows that
the curse of dimensionality is offset by the information from the
cross-sectional data when the factors are strong. Note that this result
does not need explicit constraints on the structure of
$\bSigma_{\bepsilon}$
other than implicit constraints
from assumption (D).

The assumption that all the non-zero eigenvalues of $\L$ are different is not essential,
and is merely introduced to simplify the presentation in the sense that Theorem~\ref{thm:A}
now can deal with  the convergence of the estimator for $\A$ directly. Otherwise a
discrepancy measure for two linear spaces has to be introduced in order to make statements
on the convergence rate of the estimator for the factor
loading space ${\cal M}(\A)$; see \cite{PY08}.

To present the rates of convergence for the covariance matrix estimator
$\hat{\bSigma}_{\y}$ of $\bSigma_{\y}\equiv \var(\y_t)$  and its inverse,
we introduce more conditions.

\begin{itemize}
  \item[(M1)] The error-variance matrix $\bSigma_{\bepsilon}$ is of the
form
$$ \bSigma_{\bepsilon} = \diag(\sigma_1^2\1_{m_1}^T, \sigma_2^2\1_{m_2}^T, \cdots, \sigma_k^2\1_{m_k}^T), $$
where $k, m_1, \cdots, m_k \geq 1$ are  integers,
$\1_{m_j}$ denotes the $m_j\times 1$  vector of ones,
and all the  $\sigma_j^2$ are uniformly
bounded away from 0 and infinity as $ n \rightarrow \infty $.

\item[(M2)] It holds that $p^{1-\delta}s^{-1}h_n^2 \rightarrow 0$ and $r^2s^{-1} < p^{1-\delta}h_n^2$,
where $s = \min_{1 \leq i \leq k} m_i$ and $h_n$ given in
Theorem \ref{thm:A}.
\end{itemize}
Condition (M1) assumes  that the white noise components
for $p$ time series are uncorrelated  with each other at any fixed time.
Furthermore, there are only maximum $k$ different values among their variances.
This facilitates a consistent pooled
estimator for $\bSigma_{\bepsilon}$; see also (M2).

We estimate $\bSigma_{\y}$ by
\begin{equation}\label{eqn:Sigmayhat}
\begin{split}
  \hat{\bSigma}_{\y} &= \hat{\A}\hat{\bSigma}_{\x}\hat{\A}^T + \hat{\bSigma}_{\bepsilon}, \;\; \text{ with } \;\; \hat{\bSigma}_{\x} = \hat{\A}^T(\tilde{\bSigma}_{\y} - \hat{\bSigma}_{\bepsilon})\hat{\A} \;\; \text{ and }\\
\hat{\bSigma}_{\bepsilon} &= \diag(\hat{\sigma}_1^2\1_{m_1}^T, \hat{\sigma}_2^2\1_{m_2}^T, \cdots, \hat{\sigma}_k^2\1_{m_k}^T), \;\; \hat{\sigma}_j^2 = n^{-1}s_j^{-1}\norm{\Delta_j\hat{\E}}_F,
\end{split}
\end{equation}
where $\Delta_j = \diag(\0_{m_1}^T,\cdots,\0_{m_{j-1}}^T,
\1_{m_j}^T,\0_{m_{j+1}}^T,\cdots,\0_{m_k}^T)$, $\hat{\E} = (\I_p -
\hat{\A}\hat{\A}^T)(\y_1 \cdots \y_n)$, and the norm $\norm{\M}_F$ denotes
the Frobenius norm, defined by $\norm{\M}_F = \tr(\M^T \M)^{1/2}$.

In practice we do not know the value of $k$ and the grouping.
We may start with one single group, and
estimate $\hat{\bSigma}_{\bepsilon} = \hat{\sigma}^2\I_p$. By looking at
the sample covariance matrix of the resulting  residuals,
we may group together the variables with similar magnitude of variances.
We then fit the model again with the constrained covariance
structure specified in (M1).

The theorem below presents the convergence rates for
the sample covariance estimator $\tilde{\bSigma}_{\y} \equiv \tilde{\bSigma}_{\y}(0)$ defined in
(\ref{eqn:Ltilde}) and the factor model based estimator $\hat{\bSigma}_{\y}$ defined
in (\ref{eqn:Sigmayhat}).

\begin{theorem}\label{thm:Cov}
Under assumption (D), it holds that
\begin{align*}
\norm{\tilde{\bSigma}_{\y} - \bSigma_{\y}} &= O_P(p^{1-\delta}h_n) =
O_P(p^{1-\delta}n^{-l_x} + p^{1-\delta/2}n^{-l_{x\epsilon}} + pn^{-l_{\epsilon}}).
\end{align*}
Furthermore,
$$ \norm{\hat{\bSigma}_{\y} - \bSigma_{\y}} = O_P(p^{1-\delta}h_n), $$
provided that
the condition of Theorem \ref{thm:A}, and
(M1) and (M2) also hold.
\end{theorem}

Theorem~\ref{thm:Cov} indicates that asymptotically
there is little difference in
using the sample covariance matrix or the factor model-based
covariance matrix estimator  even when all the
factors are strong, as  both the estimators have rates of
convergence linear in $p$. This result is in line with \cite{FFL08}
which shows  that the sample covariance matrix as well as the
factor model-based covariance matrix estimator are consistent in
Frobenius norm at a rate linear in $p$, with the factors known in
advance. We further illustrate this phenomenon numerically in section
\ref{sect:simulations}.

However as for the estimation for
the precision matrix $\bSigma_{\y}^{-1}$,
the estimator $\hat{\bSigma}_{\y}^{-1}$ performs significantly better than
the sample counterpart.

\begin{theorem}\label{thm:Invcov}
Under the condition of Theorem \ref{thm:A}, (M1) and (M2),
it holds  that
\begin{align*}
\norm{\hat{\bSigma}_{\y}^{-1} - \bSigma_{\y}^{-1}} &=
O_P((1+(p^{1-\delta}s^{-1})^{1/2})h_n)\\ &=
O_P(h_{n}+(ps^{-1})^{1/2}(p^{-\delta/2}n^{-l_x} + n^{-l_{x\epsilon}}
+ p^{\delta/2}n^{-l_{\epsilon}})).
\end{align*}
Furthermore,
$$ \norm{\tilde{\bSigma}_{\y}^{-1} - \bSigma_{\y}^{-1}} =
O_P(p^{1-\delta}h_n) = \norm{\tilde{\bSigma}_{\y} - \bSigma_{\y}}
$$ provided $p^{1-\delta}h_n = o(1)$ and $p < n$.
\end{theorem}

Note that if $p \geq n$, $\tilde{\bSigma}_{\y}$ is singular and the rate for the inverse sample covariance matrix becomes unbounded.
If the factors are weak
{ (i.e. $\delta>0$)},
$p$ will still be in
the above rate for $\hat{\bSigma}_{\y}^{-1}$.
{
On the other hand if the factors are strong (i.e. $\delta = 0$) and
$s \asymp p$, the rate
is independent of $p$. Hence the factor model-based estimator for the precision matrix
$\bSigma_y^{-1}$
is consistent in spectral norm irrespective of the
dimension of the problem.
 Note that the condition $s \asymp p$ is fulfilled when the number of
groups with different error variances is small.
In this case, the above  rate is better than the Frobenius norm convergence rate obtained in
Theorem 3 of \cite{FFL08}. This is not surprising since
the spectral norm is always smaller than the Frobenius norm. On the
other hand, the sample precision matrix has the convergence rate linear in
$p$, which is the same as for the sample covariance matrix.}

{
\cite{FFL08} studied the rate of convergence for a factor model-based
precision matrix estimator under the Frobenius norm and the transformed
Frobenius norm
$\norm{\cdot}_{\bSigma}$ defined as
$$ \norm{\B}_{\bSigma} = p^{-1/2}\norm{\bSigma^{-1/2}\B\bSigma^{-1/2}}_F, $$
where $\norm{\cdot}_F$ denotes the Frobenius norm. They show that the
convergence rate under the
Frobenius norm still depends on $p$,  while the rate under the
transformed Frobenius norm
$\bSigma = \bSigma_{\y}$ is independent
of $p$. Since $\bSigma_{\y}$ is unknown in practice, the latter result
has little practical impact.

\section{Factors with different levels of strength}
\label{sect:sswfactors}

\subsection{Models and two estimation procedures}

Theorems \ref{thm:A} and \ref{thm:Invcov} show that the strength of factors plays an important
role in the convergence rates of our estimators. To investigate the impact
from the presence of the different levels of factor strength, we consider the case that
the factors are of two levels of strength. The cases with more than two levels may be
treated with more complex technical details.

In view of model (\ref{eqn:factormodelWLOG}), we assume that we have two group of factors,
$\x_t = (\x_{1t}^T \; \x_{2t}^T)^T$ and $\A = (\A_1 \; \A_2)$, where $\x_{jt}$ is a $r_j \times 1$ vector and $\A_j$ is a $p \times r_j$ constant matrix for $j=1,2$. The model we consider is then
\begin{equation}\label{eqn:statfactormodelWLOG}
\begin{split}
\y_t &= \A_1\x_{1t} + \A_2\x_{2t} + \bepsilon_t, \;\; \A_j^T\A_j = \I_{r_j} \text{ and } \A_1^T\A_2 = \0, \text{ with }\\
\norm{\bSigma_{jj}(k)} &\asymp p^{1-\delta_j} \asymp \norm{\bSigma_{jj}(k)}_{\min}, \;\; \norm{\bSigma_{12}(k)} = O(p^{1-\delta_1/2 - \delta_2/2}) = \norm{\bSigma_{21}(k)},\\
\norm{\bSigma_{j\bepsilon}(k)} &= O(p^{1-\delta_j/2}),\\
\cov(\x_{js}, \bepsilon_t) &= \0 \text{ for all } s \leq t,\\
\cov(\bepsilon_s, \bepsilon_t) &= \0 \text{ for all } s\neq t.
\end{split}
\end{equation}

Unless specified otherwise, all $\y_t, \; \x_{it}$ and $\bepsilon_t$ in the sequel of this
section are defined in (\ref{eqn:statfactormodelWLOG}).

We may continue apply the estimation method outlined in section \ref{subsect:estimation} to obtain the
estimator $\hat \A \equiv ( \hat \A_1, \hat \A_2)$.
However such a simple procedure may encounter problems when some factors are weak, or are much
weaker than the others. Since the eigenvalues corresponding to those weak factors
are typically small, it may be difficult in practice to distinguish them from 0 in the presence of some
large eigenvalues.  Under those circumstances, we should remove the strong (or stronger) factors first,
and then repeat the estimation procedure again in order to identify the weak (or weaker) factors.
This is the essential idea behind the two-step procedure proposed by \cite{PP06} with some illustrative
numerical examples. We provide below a theoretical justification for using the two-step estimation
method for model (\ref{eqn:statfactormodelWLOG}) in which the factors are of different levels
of strength.

We assume that $r_1$ and $r_2$ are known. Our two-step procedure is defined as follows:
(i)~By ignoring the term $\A_2 \x_{2t}$ in model (\ref{eqn:statfactormodelWLOG}), apply
the estimation method in section~\ref{subsect:estimation} to obtain $\hat\A_1$. (ii)  By removing factor $\x_{1t}$,
\begin{equation}\label{eqn:ystar}
\y_t^* = \y_t - \hat{\A}_1\hat{\A}_1^T\y_t,
\end{equation}
estimate $\A_2$ from model $\y_t^* = \A_2 \x_{2t} + \bepsilon_t^*$ using
the method of section~\ref{subsect:estimation}.
The estimator obtained is denoted as
$\check{\A}_2$, and we denote $\check\A = (\hat\A_1, \check\A_2)$.

\subsection{Asymptotic theory}\label{subsect:theories}
In view of model (\ref{eqn:factormodelWLOG}) and the results from Lemma \ref{lemma:rates} derived from assumption (D), we directly assume the following for model (\ref{eqn:statfactormodelWLOG}):
\begin{itemize}
  \item[(D)'] For $0 \leq k \leq k_0$ and $i=1,2$, the elementwise convergence rates for
  $\tilde{\bSigma}_{ii}(k) - \bSigma_{ii}(k)$, $\tilde{\bSigma}_{12}(k) -
\bSigma_{12}(k)$, $\tilde{\bSigma}_{i\bepsilon}(k) -
\bSigma_{i\bepsilon}(k)$ and $\tilde{\bSigma}_{\bepsilon}(k) -
\bSigma_{\bepsilon}(k)$ are, respectively, $O_P(p^{1-\delta_i}n^{-l_i})$,
$O_P(p^{1-\delta_1/2-\delta_2/2}n^{-l_{12}})$, $O_P(p^{(1-\delta_i)/2}n^{-l_{i\epsilon}})$ and
$O_P(n^{-l_{\epsilon}})$, where $0 \leq l_{i},\; l_{12},\;
l_{i\epsilon},\; l_{\epsilon} \leq 1/2$. Furthermore,
$\tilde{\bSigma}_{\bepsilon i}(k) = O_P(p^{(1-\delta_i)/2}n^{-l_{i\epsilon}})$, $\tilde\bSigma_{21}(k) - \bSigma_{21}(k) = O_P(p^{1-\delta_1/2-\delta_2/2}n^{-l_{12}})$ elementwisely.
\end{itemize}

\begin{theorem}\label{thm:stationaryfactors}
Let $\norm{\bSigma_{i\bepsilon}(k)} = o(p^{1-\delta_i})$ for $i=1,2$, and condition (D)' hold.
Under model (\ref{eqn:statfactormodelWLOG}),
$$\norm{\hat{\A}_1 - \A_{1}} = O_P(\omega_1), \;\;\; \norm{\hat{\A}_2 -
\A_{2}} = O_P(\omega_2) = \norm{\hat{\A} - \A}   $$
provided $\omega_1=o(1)$ and $\omega_2=o(1)$,
and
$$ \norm{\check{\A}_2 - \A_{2}} = O_P(p^{\delta_2-\delta_1}\omega_1) = \norm{\check{\A} - \A},$$
where
\begin{align*}
  \omega_1 &=  n^{-l_{1}} + p^{\delta_1 - \delta_2}n^{-l_{2}} + p^{\frac{\delta_1-\delta_2}{2}}n^{- l_{12}}
+ p^{\delta_1/2}n^{- l_{1\epsilon}} + p^{\delta_1 - \delta_2/2}n^{- l_{2\epsilon}}
+ p^{\delta_1}n^{- l_{\epsilon}},\\
\omega_2 &= \left\{
              \begin{array}{ll}
                p^{2\delta_2 - 2\delta_1}\omega_1, & \hbox{if $\norm{\bSigma_{21}(k)}_{\min} = O(p^{1-\delta_2})$;} \\
                p^{c - 2\delta_1}\omega_1, & \hbox{if $\norm{\bSigma_{21}(k)}_{\min} \asymp p^{1-c/2}$, with $ \delta_1 + \delta_2 < c < 2\delta_2$;} \\
                p^{\delta_2 - \delta_1}\omega_1, & \hbox{if $\norm{\bSigma_{21}(k)} \asymp p^{1-\delta_1/2-\delta_2/2} \asymp \norm{\bSigma_{21}(k)}_{\min}$.}
              \end{array}
            \right.
\end{align*}
\end{theorem}

Theorem \ref{thm:stationaryfactors} indicates that while the estimators for the loading $\A_1$ on the
stronger factor $\x_{1t}$ using the two methods are exactly the same, the estimation for the loading $\A_2$ on the weaker factor may benefit from the two-step
procedure, as the convergence rate for $\check \A_2$ is faster than that for $\hat \A_2$ when
$\norm{\bSigma_{21}(k)}_{\min} = O(p^{1-\delta_2})$ for example.
The practical implication of this result is that
we should search in the residuals, after an initial fitting, for  possible weak factors,
especially if the number of non-zero eigenvalues is determined by some `eyeball' test
which remains as one of the most frequently used methods in practice.

\begin{theorem}\label{thm:comparefactormagnitude}
Under the condition of Theorem \ref{thm:stationaryfactors},
$ \norm{\hat{\x}_{2t}} = o_P(\norm{\hat{\x}_{1t}})$
 provided $p^{\delta_2 - \delta_1}\omega_1 = o(1)$, where
$\hat \x_t = ( \hat \x_{1t}^T, \; \hat \x_{2t}^T)^T$ is defined in
(\ref{eqn:estfactorresiduals}).
\end{theorem}

Theorem~\ref{thm:comparefactormagnitude} indicates that under the normalization condition
$\A \A^T= \I_r$, the different levels of factor strength will also be reflected on the magnitude
of the norms of the estimated factors. In the case that the norms of the estimated
factors, derived from the method in section~\ref{subsect:estimation},  differ
substantially, the two-step procedure may be applied to improve the
estimation.


The next two theorems are on the convergence rates for the estimation for the
covariance and the precision matrices of $\y_t$. To this end,
we recast condition (M2) for model (\ref{eqn:statfactormodelWLOG}) first.

\begin{itemize}
\item[(M2)']  For $s = \min_{1 \leq i \leq k}
m_i$, where $m_i$ are given in condition (M1), it holds that
    $$s^{-1} < p^{1-\delta_1} \norm{\A^\star - \A}^2 \; \text{ and } \;
(1 + (p^{1-\delta_1}s^{-1})^{1/2})\norm{\A^\star - \A} = o_P(1),$$
\end{itemize}
where $\A^\star$ is either $\hat \A$ or $\check \A$.

\begin{theorem}\label{thm:statcov}
Under model (\ref{eqn:statfactormodelWLOG}) and condition (D)',
\begin{align*}
\norm{\tilde{\bSigma}_{\y} - \bSigma_{\y}} &= O_P(p^{1-\delta_1}n^{-l_1}
+ p^{1-\delta_2}n^{-l_2} + p^{1-\delta_1/2}n^{-l_{1\epsilon}} +
p^{1-\delta_2/2}n^{-l_{2\epsilon}} + pn^{-l_{\epsilon}})\\
&=O_P(p^{1-\delta_1}\omega_1),\\
\norm{\hat{\bSigma}_{\y} - \bSigma_{\y}} &= O_P(p^{1-\delta_1}\norm{\hat{\A} - \A}),
\quad
\norm{\check{\bSigma}_{\y} - \bSigma_{\y}} = O_P(p^{1-\delta_1}\norm{\check{\A} - \A}),
\end{align*}
where $\tilde{\bSigma}_{\y}$ is the sample covariance matrix of $\y_t$,
 $\hat{\bSigma}_{\y}$ is the factor-model based estimator
defined by
(\ref{eqn:Sigmayhat}),
and $\check{\bSigma}_{\y}$ is defined in the same manner as $\hat{\bSigma}_{\y}$ with
$\hat{\A}$ replaced by $\check{\A}$.
\end{theorem}

Similar to Theorem \ref{thm:Cov}, the above theorem indicates that
the factor model-based approach cannot improve the estimation for the covariance
matrix of $\y_t$ over the simple sample covariance matrix. In fact, it may do worse when
the levels of factor strength differ substantially,
rendering a worse convergence rate for $\norm{\hat{\A} - \A}$.

\begin{theorem}\label{thm:statinvcov}
Let conditions (D)', (M1) and (M2)' holds.
Under model (\ref{eqn:statfactormodelWLOG}),
$$ \norm{\tilde{\bSigma}_{\y}^{-1} - \bSigma_{\y}^{-1}} = O_P(
\norm{\tilde{\bSigma}_{\y} - \bSigma_{\y}} ) $$
provided $\norm{\tilde{\bSigma}_{\y} - \bSigma_{\y}} = o_P(1)$, and
$$ \norm{\hat{\bSigma}_{\y}^{-1} -
\bSigma_{\y}^{-1}} = O_P((1+(p^{1-\delta_1} s^{-1})^{1/2})\norm{\hat{\A} - \A}), $$
$$ \norm{\check{\bSigma}_{\y}^{-1} -
\bSigma_{\y}^{-1}} = O_P((1+(p^{1-\delta_1} s^{-1})^{1/2})\norm{\check{\A} - \A}), $$
where $\tilde{\bSigma}_{\y}, \; \hat{\bSigma}_{\y}$ and $\check{\bSigma}_{\y}$ are the same
as in Theorem~\ref{thm:statcov}.
\end{theorem}


Similar to Theorem \ref{thm:Invcov}, Theorem \ref{thm:statinvcov} shows that
the factor model-based methods may improve the estimation for the inverse covariance matrix
of $\y_t$, especially with the two-step estimation method as the factors in model
(\ref{eqn:statfactormodelWLOG}) are of different levels.

\section{Simulations}\label{sect:simulations}
\setcounter{equation}{0} In this section, we illustrate our
estimation methods and their properties via two simulated examples.


{\bf Example 1}.
We start by a simple example to illustrate the properties
exhibited  in Theorem \ref{thm:A}, \ref{thm:Cov} and \ref{thm:Invcov}. Assume a one factor model
$$\y_t = \A x_t + \bepsilon_t, \;\; \epsilon_{tj} \sim \text{ i.i.d. } N(0,2^2),$$
where the factor loading $\A$ is a $p\times 1$ vector with $2\cos(2\pi i/p)$ as its $i$-th element,
 and the factor time series is defined as $x_t = 0.9x_{t-1} + \eta_t$,  where $\eta_t$  are independent
$N(0,2^2)$ random variables.
Hence we have a strong factor for this model with $\delta =0$. We set $n=200, 500$ and $p=
20, 180, 400, 1000$.  For each $(n, p)$ combination,
 we generate from the model 50 samples  and calculate the estimation
errors as in Theorems \ref{thm:A}, \ref{thm:Cov} and \ref{thm:Invcov}.
The results with $n=200$ are listed in Table~\ref{table:toyexample} below.
The results with $n=500$ is similar and thus not displayed.

\begin{table}[htbp]
\centering
  \begin{tabular}{c|c|cccc}
  \hline
   $n=200$      & $\norm{\hat{\A} - \A}$  &
$\norm{\tilde{\bSigma}_{\y}^{-1} - \bSigma_\y^{-1}}$     &  $
\norm{\hat{\bSigma}_\y^{-1} - \bSigma_\y^{-1}} $    & $ \norm{\tilde{\bSigma}_{\y} -
\bSigma_\y} $  & $ \norm{\hat{\bSigma}_\y - \bSigma_\y} $  \\
             \hline
   $p=20$   & $.022_{(.005)}$   &   $.24_{(.03)}$    &  $.009_{(.002)}$  & $218_{(165)}$  &  $218_{(165)}$ \\
  $p=180$ & $.023_{(.004)}$  &   $79.8_{(29.8)}$    &  $.007_{(.001)}$ & $1962_{(1500)}$   &  $1963_{(1500)}$ \\
  $p=400$    & $.022_{(.004)}$  &   -           &  $.007_{(.001)}$  & $4102_{(3472)}$   &  $4103_{(3471)}$  \\
  $p=1000$ &   $.023_{(.004)}$  &   -          &   $.007_{(.001)}$    & $10797_{(6820)}$     & $10800_{(6818)}$       \\
  \hline
\end{tabular}
\caption{\emph{Mean of estimation errors for the one factor example.
The numbers in brackets are the corresponding standard deviations.}}
\label{table:toyexample}
\end{table}

It is
clear from Table \ref{table:toyexample} that the estimation errors
in $L_2$ norm for $\hat{\A}$ and $\hat{\bSigma}_\y^{-1}$ are
independent of $p$, as indicated by Theorem
\ref{thm:A} and Theorem \ref{thm:Invcov} with $\delta = 0$, as we
have a strong factor in this example. The inverse of the sample
covariance matrix $\tilde{\bSigma}_{\y}^{-1}$ is not defined for $p > n$, and is a bad
estimator even when $p < n$ as seen in above table. The last two
columns show that the errors of estimators $\tilde{\bSigma}_{\y}$ and $\hat{\bSigma}_\y$ increase
as $p$ increases. This is  in agreement with
Theorem \ref{thm:Cov}.

\vspace{12pt}

{\bf Example 2}.
Now we consider a model with factors of different levels of strength.
We generate data from model
(\ref{eqn:statfactormodelWLOG}) with $r=3$ factors:
\begin{align*}
  x_{1,t} &= -0.8x_{1,t-1} + 0.9e_{1,t-1} + e_{1,t},\\
  x_{2,t} &= -0.7x_{2,t-1} + 0.85e_{2,t-1} + e_{2,t},\\
  x_{3,t} &= \;\;\: 0.8x_{2,t} - 0.5x_{3,t-1} + e_{3,t},
\end{align*}
where $e_{i,t}$ are independent $ N(0,1)$ random variables.
For each column of $\A$, we generate the first $p/2$ elements
randomly from the $U(-2,2)$ distribution; the rest are set to zero.
We then adjust the strength of the
factors by normalizing the columns, setting $\a_i/p^{\delta_i/2}$ as
the $i$-th column of $\A$ (we set $\delta_2 = \delta_3$).
We let $\bepsilon_{t}$ be $p\times 1$ independent random vectors with
mean 0 and variance diag$(0.5, 0.8, 0.5, 0.8, \cdots )$, and the distributions
of all the components of $\bepsilon_{t}$ are either normal or $t_5$ (properly
normalized such that the variance is either 0.5 or 0.8).

We set $n=100, 200, 500, 1000$ and $p=100, 200, 500$.
The first factor has strength index $\delta_1$ and the last two factors
have strength index $\delta_2$. Both  $\delta_1$ and $\delta_2$ take
values $0,0.5$ or 1.
For each combination of $(n, p, \delta_1, \delta_2)$, we replicate
the simulation 100 times, and calculate the mean and the standard deviations
of the error measures.


\begin{figure}[ht]
  \includegraphics[width=1.1\textwidth, height=0.9\textwidth]{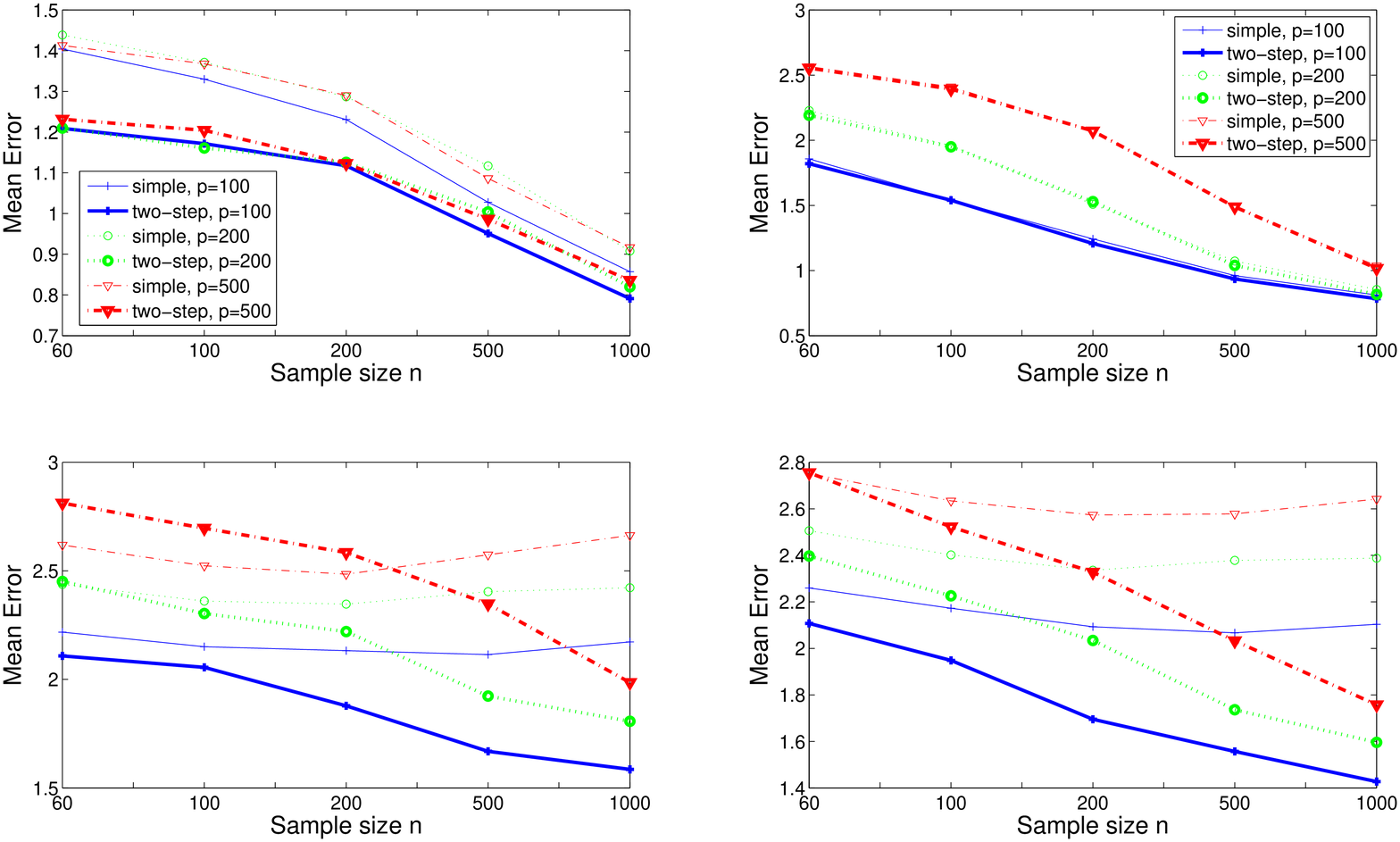}
  \caption{\emph{Mean of $\norm{\hat{\bSigma}_{\y}^{-1} -
\bSigma_{\y}^{-1}}$ and $\norm{\check{\bSigma}_{\y}^{-1} -
\bSigma_{\y}^{-1}}$, noises distributed as $t_5$. Top row:
$\delta_1=\delta_2=\delta_3=0$.
  Bottom row: $\delta_1=0, \delta_2=\delta_3=0.5$. Left column: Three factors used in estimation. Right column: Four factors used in estimation. The legends for the bottom row are the same as the top row.}}\label{figure:1}
\end{figure}

Figure \ref{figure:1} displays the mean of $\norm{\hat{\bSigma}_{\y}^{-1} - \bSigma_{\y}^{-1}}$
in the 100 replications.
For precision matrix estimation, figure \ref{figure:1} shows clearly that
when the number of factors are not underestimated, the two step procedure
outperforms the simple one when strength of factors are different, and
performs at least as good when the factors are of the same strength. The
performance is better when the number of factors is in fact more than the
optimal because we have used $k_0=3$ instead of just 1 or 2 when the
serial correlations for the factors are in fact quite weak. Hence we
accumulate pure noises, which sometimes introduces non-genuine factors
that are stronger than the genuine ones, and requires the inclusion of
more than necessary factors to reduce the errors. Not shown here, we have
repeated the simulations with $k_0=1$, and the performance is much better
and is optimal when the number of factors used is 3. The two-step procedure
still outperforms the simple one. The simulations with normal errors are
not shown here since the results are similar. The mean of $\norm{\hat\A - \A}$ and $\norm{\check\A - \A}$ exhibit similar patterns as shown in figure \ref{figure:1} for $\norm{\hat{\bSigma}_{\y}^{-1} -
\bSigma_{\y}^{-1}}$ and $\norm{\check{\bSigma}_{\y}^{-1} -
\bSigma_{\y}^{-1}}$ respectively, and the results are not shown.

For the covariance matrix estimation, our results (not shown) show
that, as in Theorem \ref{thm:statcov}, both the sample covariance
matrix and factor model based one are poor estimators when $p$ is
large. In fact both the simple and two-step procedures yield worse
estimation errors than the sample covariance matrix, although
performance gap closes down as $n$ gets larger.

\section{Data Analysis : Implied Volatility Surfaces}\label{sect:dataanalysis}
We illustrate the methodology developed through modeling the dynamic behavior of IBM, Microsoft and Dell implied volatility surfaces.
The data was obtained from OptionMetrics via the WRDS database. The dates in question are $03/01/2006 - 29/12/2006$ (250 days in total). For each day $t$ we observe the implied volatility $W_t(u_i, v_j)$ computed from call options as a function of time to maturity of 30, 60, 91, 122, 152, 182, 273, 365, 547 and 730 calender days which we denote by $u_i$, $i=1,\dots,p_u$ ($p_u = 10$) and deltas of 0.2, 0.25, 0.3, 0.35, 0.4, 0.45, 0.5, 0.55, 0.6, 0.65, 0.7, 0.75, and 0.8 which we denote by $v_j$, $j=1,\dots,p_v$ ($p_v = 13$). We collect these implied volatilities in the matrix $\W_t = (W_t(u_i, v_i)) \in \mathbb{R}^{p_u \times p_v}$. Figure \ref{figure:3g} displays the mean volatility surface of IBM, Microsoft and Dell over the period in question. It is clear from this graphic that the implied volatilities surfaces are not flat. Indeed any cross-section in the maturity or delta axis display the well documented volatility smile.

\begin{figure}[!htp]
  \begin{center}
      \includegraphics[width=1.6in]{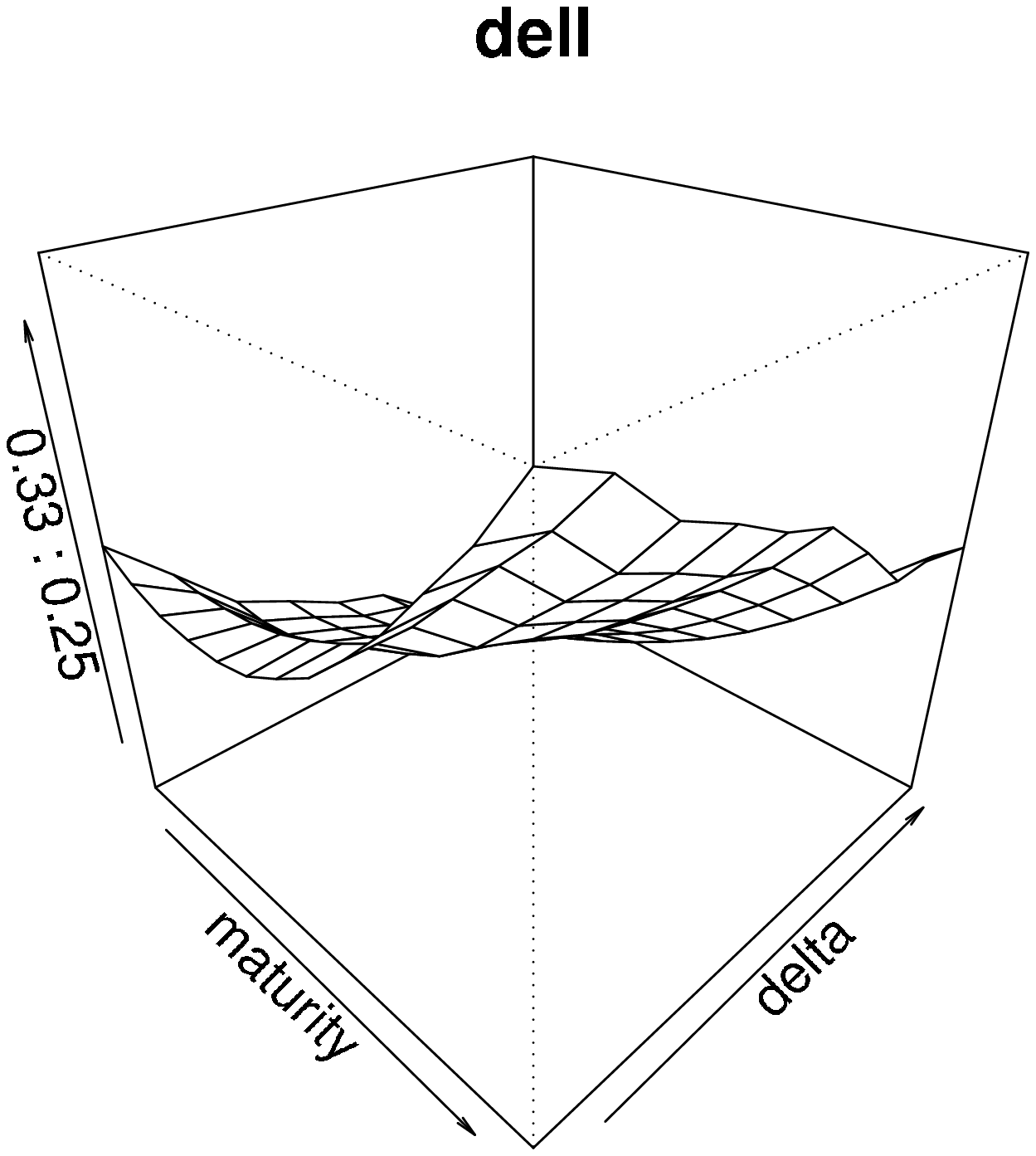}
      \includegraphics[width=1.6in]{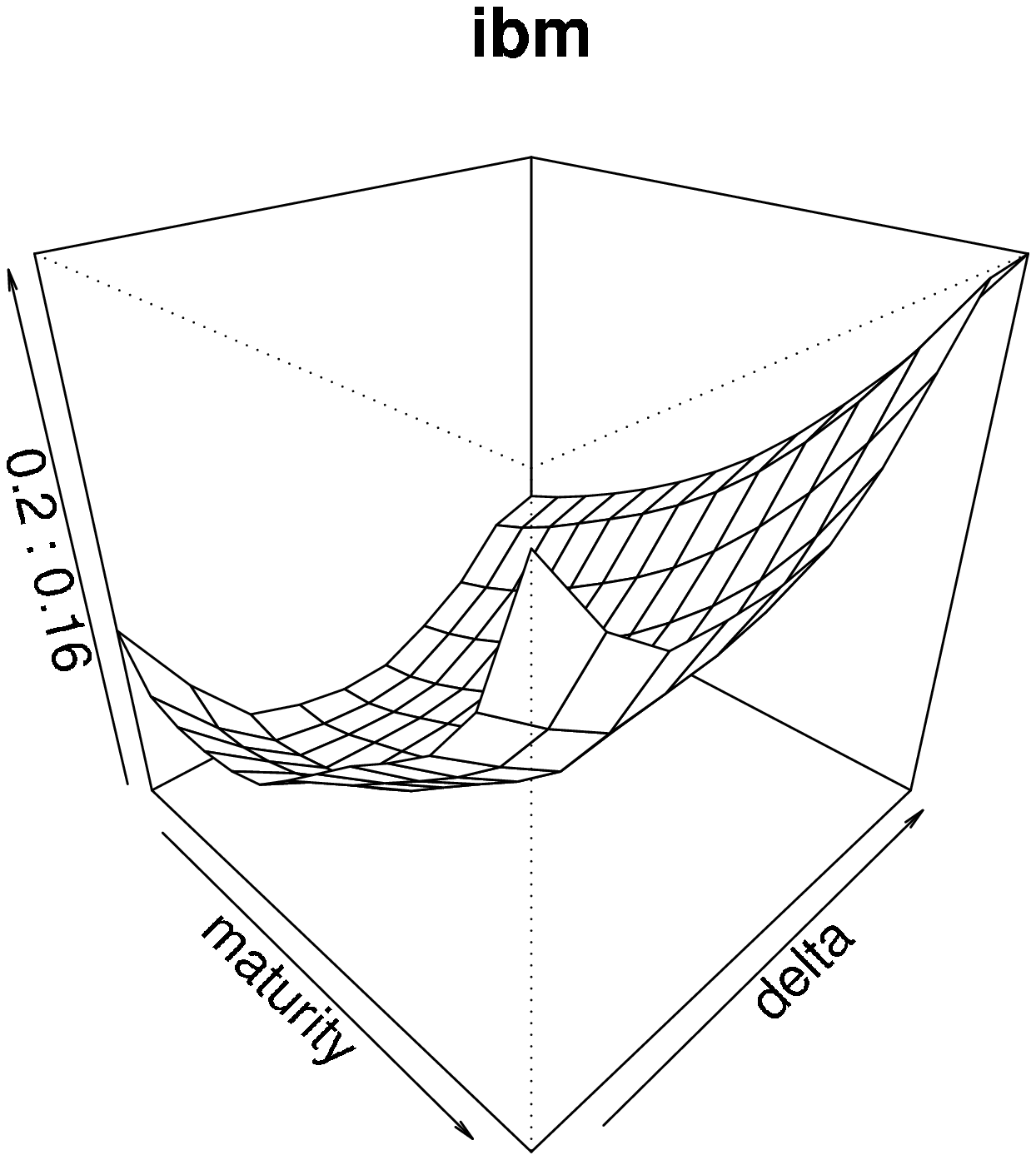}
        \includegraphics[width=1.6in]{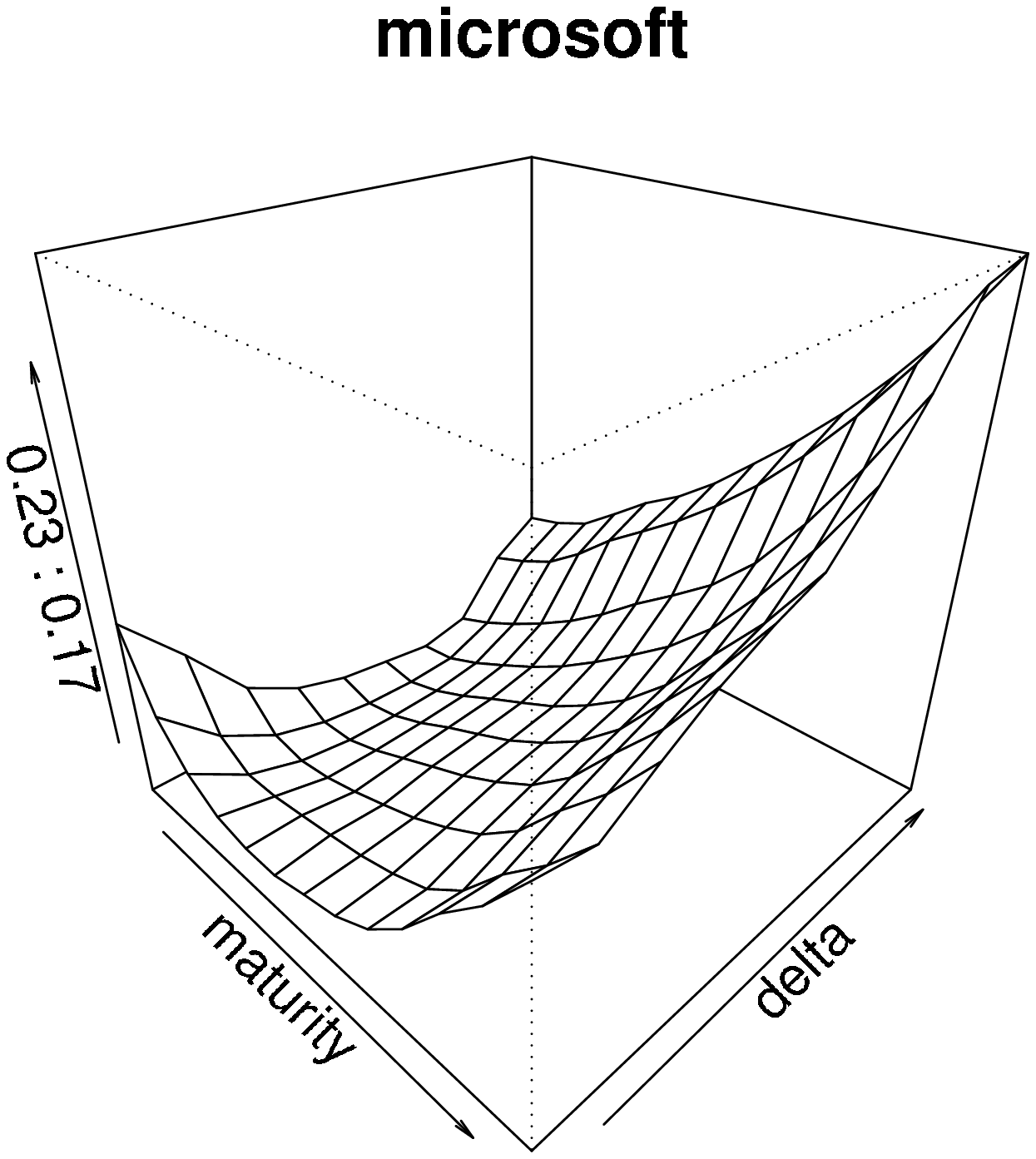}
\caption[Fig 8]{\sl Mean implied volatility surfaces.}\label{figure:3g}.
\end{center}
\end{figure}
It is a well documented stylized fact that implied volatilities are non-stationary (see  \cite{Cd02}, \cite{FHM07} and \cite{PMHB09} amongst others). Indeed, when applying the Dickey-Fuller test to each of the univariate time series $W_t(u_i, v_i)$, none of the $p_u \times p_v = 130$ nulls of unit roots could be rejected at the 10$\%$ level. 
Of course we should treat the results of these tests with some
caution since we are performing a large number of hypothesis tests,
but even still the evidence in favor of unit roots is overwhelming.
Therefore, instead of working with $\W_t$ directly, we choose to
work with $\Delta \W_t = \W_t - \W_{t-1}$. Our observations are then
$\y_t = \mbox{vec} \{ \Delta \W_t \}$, where for any matrix $\M =
(\m_1,\dots,\m_{p_v}) \in \mathbb{R}^{p_u \times p_v}$,
$\mbox{vec}\{\M \} = (\m_1^T,\dots,\m_{p_v}^T)^T \in \mathbb{R}^{p_u
p_v}$.
Note that $\y_t$ is now defined over $04/01/2006 -
29/12/2006$ since we lose an observation due to differencing. Hence altogether there are 249 time points, and the dimension of $\y_t$ is
$p = p_v \times p_u = 130$.

We perform the factor model estimation on a rolling window of length 100 days. A window is defined from the $i$-th day to the $(i + 99)$-th day for $i=1,\cdots,150$. The length of the window is chosen so that the stationary assumption of the data is approximately satisfied. For each window, we compare our methodology with the least squares based methodology by \cite{BN02} by estimating the factor loadings matrix and the factors series for the two methods. For the $i$-th window, we use an AR model to forecast the $(i+100)$-th value of the estimated factor series $\x_{i+100}^{(1)}$, so as to obtain a one-step ahead forecast $\y_{i+100}^{(1)} = \hat{\A}\x_{i+100}^{(1)}$ for $\y_{i+100}$. We then calculate the RMSE for the $(i+100)$-th day defined by
$$ \text{RMSE} = p^{-1/2} \norm{\y_{i+100}^{(1)} - \y_{i+100}}.$$
More in depth theoretical as well as data analysis for forecasting is given in \cite{LYB10b}.

\subsection{Estimation results}
In forming the matrix $\tilde\L$ for each window, we take
$k_0=5$ in (\ref{eqn:Ltilde}) , taking advantage that the autocorrelations are not weak
even at higher lags, though similar results (not reported here) are
obtained for smaller $k_0$.

\begin{figure}[!htp]
  \begin{center}
      \includegraphics[width=3.5in]{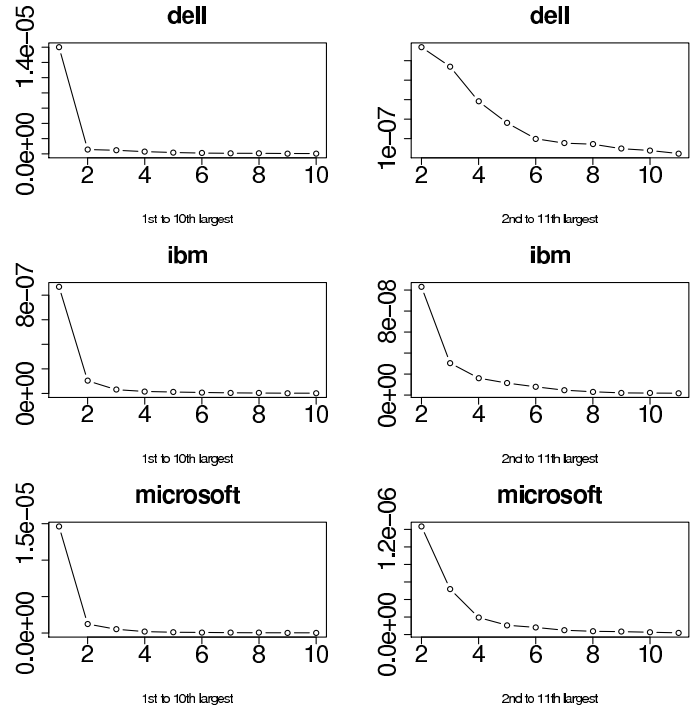}
\caption[Fig 8]{\sl Averages of ordered eigenvalues of $\tilde\L$ over the 150 windows. Left: Ten largest. Right: Second to eleventh largest.}\label{figure:3i}
\end{center}
\end{figure}
Figure \ref{figure:3i} displays the average of each ordered eigenvalue over the 150 windows. The left hand side shows the average of the
largest to the average of the tenth largest eigenvalue of $\tilde\L$ for Dell, IBM and Microsoft for our method, whereas the right hand side shows the second to eleventh largest. We obtain similar results for the \cite{BN02} procedure and thus the corresponding graph is not shown.

From this graphic it is apparent that there is one eigenvalue that is
much larger than the others for all three companies for each window. We have done automatic selection for the number of factors for each window using the $IC_{p1}$ criterion in \cite{BN02} and a one factor model is consistently obtained for each window and for each company. Hence both methods chose a one factor model over the 150 windows.

\begin{figure}[!htp]
  \begin{center}
      \includegraphics[width=2in]{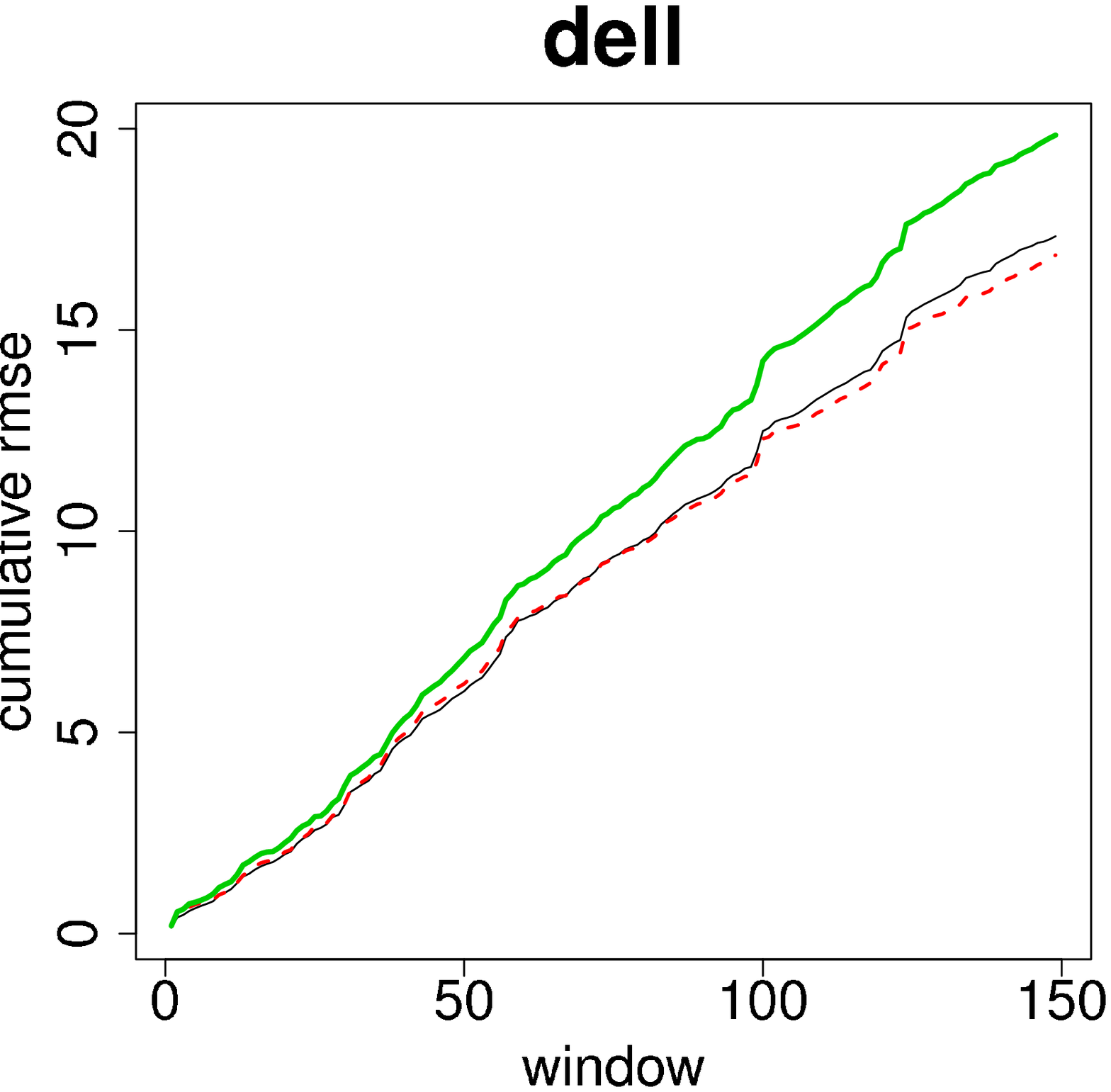}
      \includegraphics[width=2in]{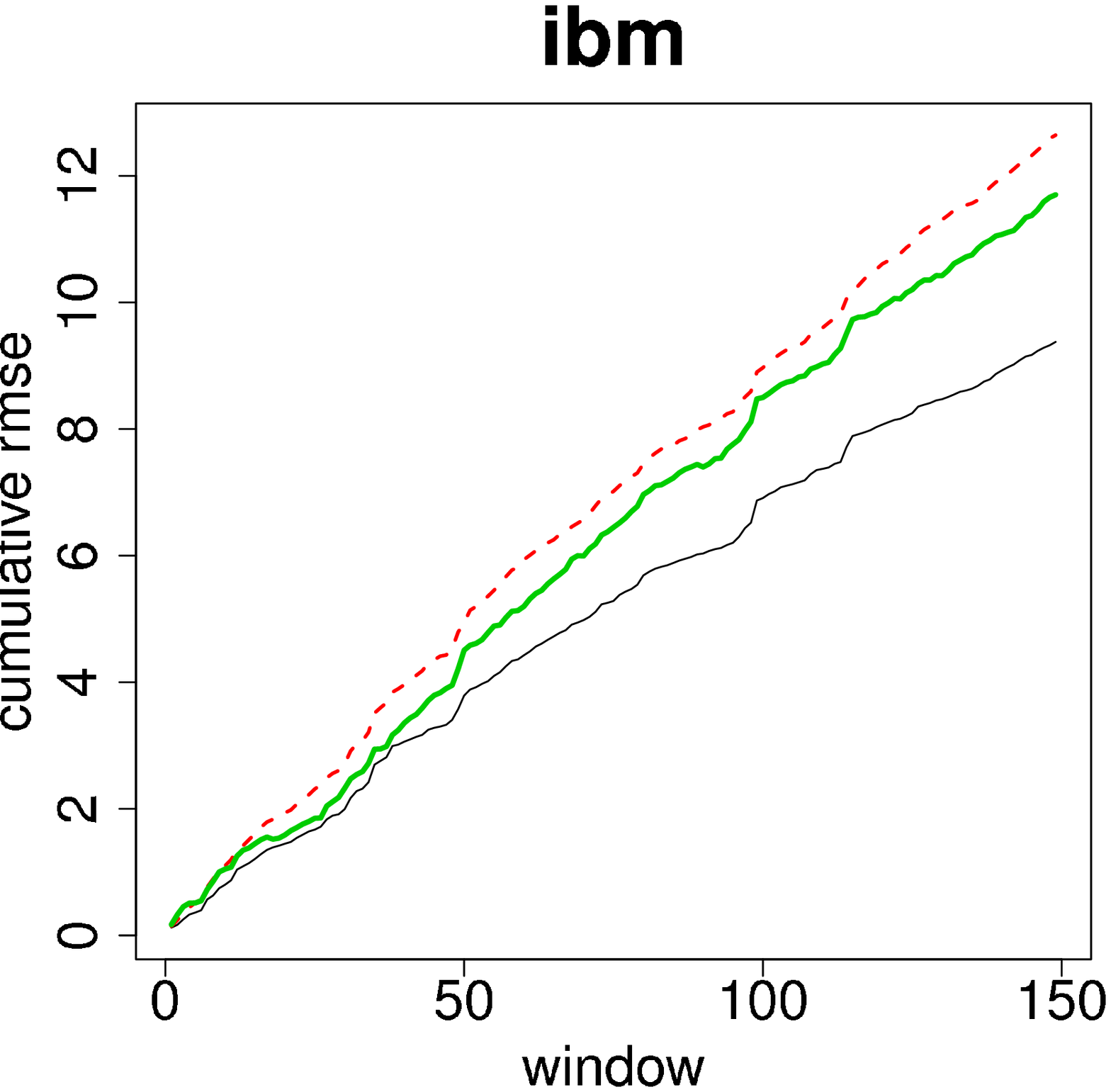}
        \includegraphics[width=2in]{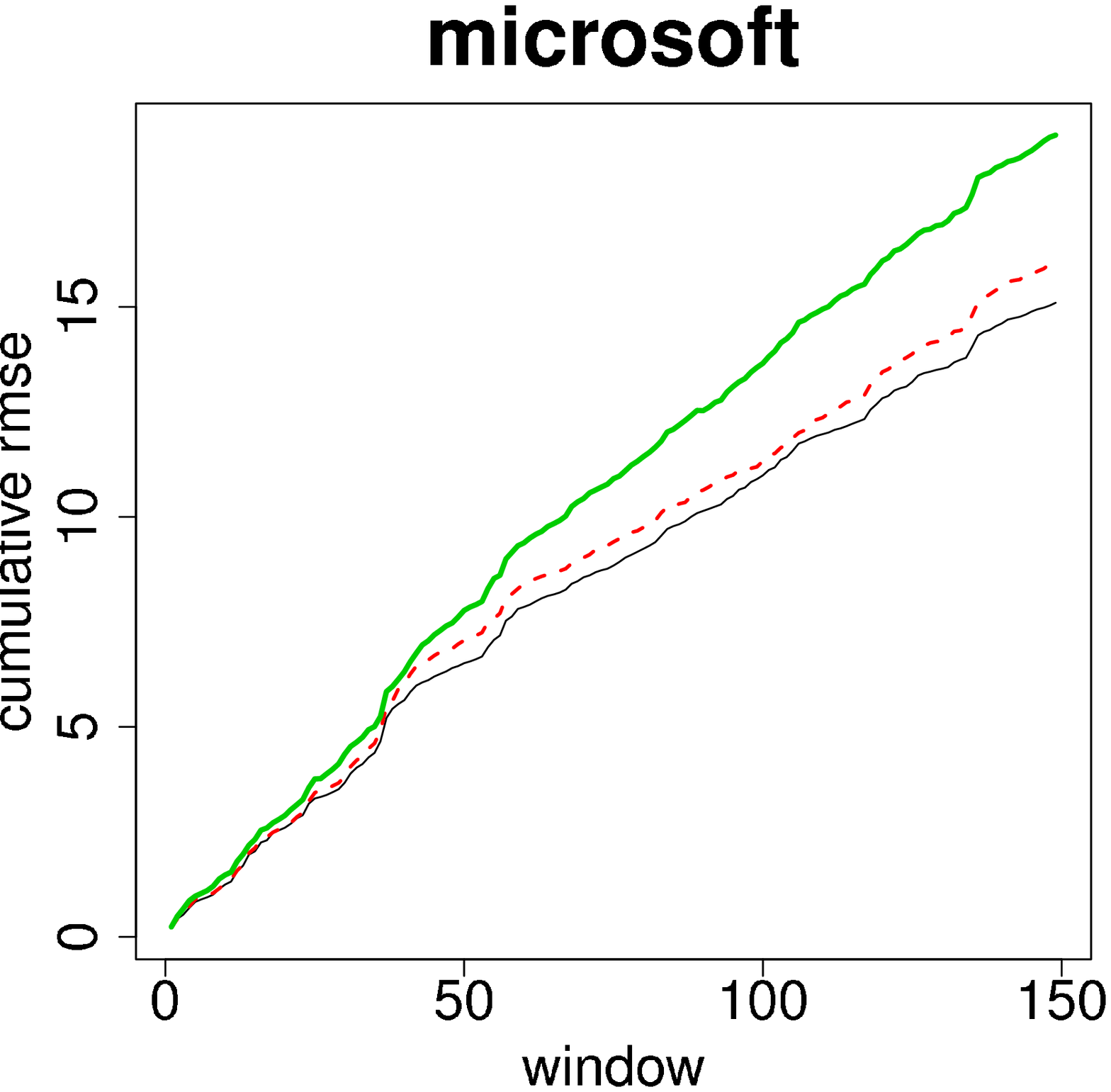}
\caption[Fig 8]{\sl The cumulative RMSE over the 150 windows. \textcolor[rgb]{0.98,0.00,0.00}{Red dotted:}  \cite{BN02} procedure. \textcolor[rgb]{0.00,1.00,0.00}{Green: } Taking forecast $\y_{t+1}^{(1)}$ to be $\y_t$. Black: Our method.}\label{figure:3j}
\end{center}
\end{figure}

Figure \ref{figure:3j} displays the cumulative RMSE over the 150 windows for each method. We choose a benchmark procedure (green line in each plot), where we just treat today's value as the one-step ahead forecast. Except for Dell where \cite{BN02} procedure is doing marginally better, our methodology consistently outperforms the benchmark procedure and is better than \cite{BN02} for IBM and Microsoft.

\section{Proofs}\label{sect:Proofs}
\setcounter{equation}{0}
First of all, we show how model (\ref{eqn:factormodelWLOG}) can be derived from (\ref{eqn:factormodel}).

Applying the standard
QR decomposition, we may write $\A = \Q\R$, where $\Q$
is a $p \times r$ matrix such that $\Q^T\Q = \I_r$, $\R$ is an $r
\times r$ upper triangular matrix. Therefore model
(\ref{eqn:factormodel}) can be expressed as
\begin{equation*}
\y_t = \Q\x_t^\prime + \bepsilon_t,
\end{equation*}
where $\x_t^\prime = \R\x_t$. With assumptions (A) to (C), the
diagonal entries of $\R$ are all asymptotic to $p^{\frac{1-\delta}{2}}$. Since $r$ is a constant, using
$$\norm{\R} = \max_{\norm{\u}=1} \norm{\R \u}, \;\;\; \norm{\R}_{\min} = \min_{\norm{\u}=1} \norm{\R \u},$$
and
the fact that $\R$ is an $r \times r$ upper triangular matrix with all diagonal elements having the largest order $p^{\frac{1-\delta}{2}}$, we have
$$\norm{\R} \asymp p^{\frac{1-\delta}{2}} \asymp \norm{\R}_{\min}.$$
Thus,
for $k = 1,\cdots, k_0$,
$\bSigma_{\x^\prime}(k) = \cov(\x_{t-k}^\prime, \x_t^\prime) = \R\bSigma_{\x}(k)\R^T ,$
with
$$ p^{1-\delta} \asymp \norm{\R}_{\min}^2 \cdot \norm{\bSigma_{\x}(k)}_{\min} \leq \norm{\bSigma_{\x^\prime}(k)}_{\min} \leq  \norm{\bSigma_{\x^\prime}(k)} \leq \norm{\R}^2 \cdot \norm{\bSigma_{\x}(k)}
\asymp p^{1-\delta}, $$
so that $\norm{\bSigma_{\x}(k)} \asymp p^{1-\delta} \asymp \norm{\bSigma_{\x}(k)}_{\min}$.
We used $\norm{\A\B}_{\min} \geq \norm{\A}_{\min}\cdot\norm{\B}_{\min}$, which can be proved by noting
\begin{align}
  \norm{\A\B}_{\min} &= \min_{\u \neq \0} \frac{\u^T\B^T\A^T\A\B\u}{\norm{\u}^2} \geq \min_{\u \neq \0} \frac{(\B\u)^T\A^T\A(\B\u)}{\norm{\B\u}^2} \cdot \frac{\norm{\B\u}^2}{\norm{\u}^2} \notag\\
  &\geq \min_{\w \neq \0}\frac{\w^T\A^T\A\w}{\norm{\w}^2} \cdot \min_{\u \neq \0} \frac{\norm{\B\u}^2}{\norm{\u}^2} = \norm{\A}_{\min}\cdot\norm{\B}_{\min}. \label{eqn:normmininq}
\end{align}
Finally, using assumption (A) that $\bSigma_{\x,\bepsilon}(k) = O(1)$ elementwise, and that it has $rp \asymp p$ elements, we have
\begin{align*}
  \norm{\bSigma_{\x^\prime,\bepsilon}(k)} = \norm{\R\bSigma_{\x,\bepsilon}(k)} \leq \norm{\R} \cdot \norm{\bSigma_{\x,\bepsilon}(k)}_F = O(p^{\frac{1-\delta}{2}}) \cdot O(p^{1/2}) = O(p^{1-\delta/2}).
\end{align*}

Before proving the theorems in section \ref{sect:theories}, we need to have three lemmas.

\begin{lemma}\label{lemma:rates}
Under the factor model (\ref{eqn:factormodelWLOG}) which is a reformulation of (\ref{eqn:factormodel}), and under condition (D) in section \ref{subsect:identifiability}, we have for $0\le k \le k_0$,
\begin{align*}
\norm{\tilde{\bSigma}_{\x}(k) - \bSigma_{\x}(k)} &= O_P(p^{1-
\delta}n^{-l_x}), \;\;\quad \norm{\tilde{\bSigma}_{\bepsilon}(k) - \bSigma_{\bepsilon}(k)} = O_P(pn^{-l_{\epsilon}}),\\
 \norm{\tilde{\bSigma}_{\x, \bepsilon}(k) - \bSigma_{\x,\bepsilon}(k)} &= O_P(p^{1-\delta/2}n^{-l_{x\epsilon}})
 = \norm{\tilde{\bSigma}_{\bepsilon, \x}(k) - \bSigma_{\bepsilon, \x}(k)},
\end{align*}
for some constants $0 < l_x, l_{x\epsilon}, l_{\epsilon} \leq 1/2$. Moreover, $\norm{\x_t}^2 = O_P(p^{1-\delta})$ for all real $t$.
\end{lemma}

{\bf Proof.} Using the notations in section \ref{subsect:identifiability}, let $\x_t$ be the factors in model (\ref{eqn:factormodel}), and $\x_t^\prime$ be the factors in model (\ref{eqn:factormodelWLOG}), with the relation that $ \x_t^\prime = \R\x_t $, where $\R$ is an upper triangular matrix with $\norm{\R} \asymp p^{\frac{1-\delta}{2}} \asymp \norm{\R}_{\min}$ (see the start of this section for more details on $\R$). Then we immediately have $\norm{\x_t^{\prime}}^2 \leq \norm{\R}^2 \cdot \norm{\x_t}^2 = O_P(p^{1-\delta}r) = O_P(p^{1-\delta})$.

Also, the covariance matrix and the sample covariance matrix for $\{\x_t^\prime\}$ are respectively
$$ \bSigma_{\x^\prime}(k) = \R\bSigma_{\x}(k)\R^T, \;\;\; \tilde{\bSigma}_{\x^\prime}(k) = \R\tilde{\bSigma}_{\x}(k)\R^T, $$
where $\bSigma_{\x}(k)$ and $\tilde{\bSigma}_{\x}(k)$ are respectively the covariance matrix and the sample covariance matrix for the factors $\{\x_t\}$. Hence
\begin{align*}
\norm{\tilde{\bSigma}_{\x^\prime}(k) - \bSigma_{\x^\prime}(k)} &\leq \norm{\R}^2 \cdot \norm{\tilde{\bSigma}_{\x}(k) - \bSigma_{\x}(k)}\\
&=O(p^{1-\delta}) \cdot O_P(n^{-l_x} \cdot r)\\
&=O_P(p^{1-\delta}n^{-l_x}),
\end{align*}
which is the rate specified in the lemma. We used the fact that the matrix $\tilde{\bSigma}_{\x}(k) - \bSigma_{\x}(k)$ has $r^2$ elements, with elementwise rate of convergence being $O(n^{-l_x})$ as in assumption (D). Other rates can be derived similarly. $\square$

The following is Theorem 8.1.10 in \cite{GV96}, which is stated explicitly since most of our main theorems are based on this. See \cite{JL09} also.

\begin{lemma}\label{lemma:matrixperturbation}
Suppose $\A$ and $\A + \E$ are $n \times n$ symmetric matrices and that
$$ \Q = [\Q_1 \;\; \Q_2] \;\;\; (\Q_1 \; \text{ is } \; n \times r, \; \Q_2 \; \text{ is } \; n \times (n-r)) $$ is an orthogonal matrix such that {\em span($\Q_1$)} is an invariant subspace for $\A$ (i.e., {\em span($\Q_1) \subset$ span($\A$)}). Partition the matrices $\Q^T\A\Q$ and $\Q^T\E\Q$ as follows:
\begin{align*}
  \Q^T\A\Q = \left(
               \begin{array}{cc}
                 \D_1 & \0 \\
                 \0 & \D_2 \\
               \end{array}
             \right) \;\;\;\;\;\; \Q^T\E\Q = \left(
               \begin{array}{cc}
                 \E_{11} & \E_{21}^T \\
                 \E_{21} & \E_{22} \\
               \end{array}
               \right).
\end{align*}
If {\em sep}$(\D_1, \D_2) := \min_{{\lambda \in \lambda(\D_1), \; \mu \in \lambda(\D_2)}}|\lambda - \mu| > 0$, where $\lambda(M)$ denotes the set of eigenvalues of the matrix $M$, and
$$ \norm{\E} \leq \frac{\text{{\em sep}}(\D_1, \D_2)}{5}, $$
then there exists a matrix $\P \in \mathbb{R}^{(n-r) \times r}$ with
$$ \norm{\P} \leq \frac{4}{\text{{\em sep}}(\D_1, \D_2)}\norm{\E_{21}} $$
such that the columns of $\hat{\Q}_1 = (\Q_1 + \Q_2\P)(\I + \P^T\P)^{-1/2}$ define an orthonormal basis for a subspace that is invariant for $\A + \E$.
\end{lemma}

{\bf Proof of Theorem \ref{thm:A}.} Under model (\ref{eqn:factormodelWLOG}), the assumption that $\norm{\bSigma_{\x,\bepsilon}(k)} = o(p^{1-\delta})$, and the definition of $\L$ and $\D_x$ in section \ref{subsect:estimation} such that $\L\A = \A\D$,\; $\D$ has non-zero eigenvalues of order $p^{2-2\delta}$, contributed by the term $\bSigma_{\x}(k)\bSigma_{\x}(k)^T$. If $\B$ is an orthogonal complement of $\A$, then $\L\B = \0$, and
\begin{align}\label{eqn:thm:A:sepDx}
\left(
     \begin{array}{c}
       \A^T \\
       \B^T \\
     \end{array}
   \right)\L (\A \;\; \B) = \left(
                                  \begin{array}{cc}
                                    \D & \0 \\
                                    \0 & \0 \\
                                  \end{array}
                                \right),
\end{align}
with $\sep(\D, \0) = \lambda_{\min}(\D) \asymp p^{2-2\delta}$ (see Lemma \ref{lemma:matrixperturbation} for the definition of the function $\sep$).

Define $\E_\L = \tilde{\L} - \L$, where $\tilde{\L}$ is defined in (\ref{eqn:Ltilde}). Then it is easy to see that
\begin{equation}\label{eqn:thm:A:E_Lbound}
\norm{\E_{\L}} \leq \sum_{k=1}^{k_0} \Big\{ \norm{\tilde{\bSigma}_{\y}(k) - \bSigma_{\y}(k)}^2 + 2\norm{\bSigma_{\y}(k)} \cdot \norm{\tilde{\bSigma}_{\y}(k) - \bSigma_{\y}(k)}  \Big\}.
\end{equation}
Suppose we can show further that
\begin{equation}\label{eqn:thm:A:E_Lorder}
\norm{\E_{\L}} = O_P(p^{2-2\delta}n^{-l_x} + p^{2-3\delta/2}n^{-l_{x\epsilon}} + p^{2-\delta}n^{-l_{\epsilon}}) = O_P(p^{2-2\delta}h_n),
\end{equation}
then since $h_n = o(1)$, we have from (\ref{eqn:thm:A:sepDx}) that
$$ \norm{\E_{\L}} = O_P(p^{2-2\delta}h_n) \leq \sep(\D, \0)/5 $$ for sufficiently large $n$. Hence we can apply Lemma \ref{lemma:matrixperturbation} to conclude that there exists a matrix $\P \in \mathbb{R}^{(p-r) \times r}$ such that
$$ \norm{\P} \leq \frac{4}{\sep(\D, \0)}\norm{(\E_{\L})_{21}} \leq \frac{4}{\sep(\D, 0)}\norm{\E_\L} = O_P(h_n), $$ and $\hat{\A} = (\A + \B\P)(\I + \P^T\P)^{-1/2}$ is an estimator for $\A$. Then we have
\begin{align*}
\norm{\hat{\A} - \A} &= \norm{(\A(\I - (\I + \P^T\P)^{1/2}) + \B\P)(\I + \P^T\P)^{-1/2}}\\
&\leq \norm{\I - (\I + \P^T\P)^{1/2}} + \norm{\P}\\
&\leq 2\norm{\P} = O_P(h_n).
\end{align*}

Hence it remains to show (\ref{eqn:thm:A:E_Lorder}). To this end, consider for $k \geq 1$,
\begin{align}\label{eqn:thm:A:Sigmaynorm}
\norm{\bSigma_{\y}(k)} = \norm{\A\bSigma_{\x}(k)\A^T + \A\bSigma_{\x,\bepsilon}(k)} \leq \norm{\bSigma_{\x}(k)} + \norm{\bSigma_{\x, \bepsilon}(k)} = O(p^{1-\delta})
\end{align}
by assumptions in model (\ref{eqn:factormodelWLOG}) and $\norm{\bSigma_{\x,\bepsilon}(k)} = o(\norm{\bSigma_{\x}(k)})$. Finally, noting $\norm{\A} = 1$,
\begin{equation}\label{eqn:thm:A:Sigmaytildenorm}
\begin{split}
\norm{\tilde{\bSigma}_{\y}(k) - \bSigma_{\y}(k)} \leq &\norm{\tilde{\bSigma}_{\x}(k) - \bSigma_{\x}(k)} + 2\norm{\tilde{\bSigma}_{\x,\bepsilon}(k) - \bSigma_{\x,\bepsilon}(k)}\\
&+\norm{\tilde{\bSigma}_{\bepsilon}(k) - \bSigma_{\bepsilon}(k)}\\
=&O_P(p^{1-\delta}n^{-l_x} + p^{1-\delta/2}n^{-l_{x\epsilon}} + pn^{-l_{\epsilon}})
\end{split}
\end{equation}
by Lemma \ref{lemma:rates}. With (\ref{eqn:thm:A:Sigmaynorm}) and (\ref{eqn:thm:A:Sigmaytildenorm}), we can conclude from (\ref{eqn:thm:A:E_Lbound}) that
$$ \norm{\E_\L} = O_P(\norm{\bSigma_{\y}(k)} \cdot \norm{\tilde{\bSigma}_{\y}(k) - \bSigma_{\y}(k)}), $$
which is exactly the order specified in (\ref{eqn:thm:A:E_Lorder}). $\square$

\vspace{12pt}

{\bf Proof of Theorem \ref{thm:Cov}.} For the sample covariance matrix $\tilde{\bSigma}_{\y}$, note that (\ref{eqn:thm:A:Sigmaytildenorm}) is applicable to the case when $k=0$, so that
$$ \norm{\tilde{\bSigma}_{\y} - \bSigma_{\y}} = O_P(p^{1-\delta}n^{-l_x} + p^{1-\delta/2}n^{-l_{x\epsilon}} + pn^{-l_{\epsilon}}) = O_P(p^{1-\delta}h_n). $$
For $\hat{\bSigma}_{\y}$, we have
\begin{equation}\label{eqn:thm:Cov:Sigmayhatnorm}
\begin{split}
\norm{\hat{\bSigma}_{\y} - \bSigma_{\y}} &\leq \norm{\hat{\A}\hat{\bSigma}_{\x}\hat{\A}^T - \A\bSigma_{\x}\A^T} + \norm{\hat{\bSigma}_{\bepsilon} - \bSigma_{\bepsilon}}\\
&= O_P(\norm{\hat{\A} - \A}\cdot\norm{\hat{\bSigma}_{\x}} + \norm{\hat{\bSigma}_{\x} - \bSigma_{\x}} + \norm{\hat{\bSigma}_{\bepsilon} - \bSigma_{\bepsilon}}).
\end{split}
\end{equation}
We first consider $\norm{\hat{\bSigma}_{\bepsilon} - \bSigma_{\bepsilon}} = \max_{1 \leq j \leq p}|\hat{\sigma}_j^2 - \sigma_j^2| \leq I_1 + I_2 + I_3$, where
\begin{align*}
 I_1 &= n^{-1}s^{-1}\norm{\Delta_j(\I_p - \hat{\A}\hat{\A}^T)\A\X}_F^2,\\
 I_2 &= |n^{-1}s^{-1}\norm{\Delta_j(\I_p - \hat{\A}\hat{\A}^T)\E}_F^2 - \sigma_j^2|,\\
 I_3 &= 2n^{-1}s^{-1}\norm{\Delta_j(\I_p - \hat{\A}\hat{\A}^T)\A\X}_F \cdot \norm{\Delta_j(\I_p - \hat{\A}\hat{\A}^T)\E}_F,
\end{align*}
with $\X = (\x_1 \cdots \x_n)$, $\E = (\bepsilon_1 \cdots \bepsilon_n)$. We have
\begin{align}
  I_1 &\leq n^{-1}s^{-1}\norm{\Delta_j(\I_p - \hat{\A}\hat{\A}^T)(\A - \hat{\A})\X}_F^2 \notag\\
      &\leq n^{-1}s^{-1} \norm{\Delta_j}^2 \cdot \norm{\I_p - \hat{\A}\hat{\A}^T}^2 \cdot \norm{\A - \hat{\A}}^2 \cdot \norm{\X}_F^2 \notag\\
      &=O_P(p^{1-\delta}s^{-1}h_n^2), \label{eqn:thm:Cov:I1}
\end{align}
where we used Theorem \ref{thm:A} for the rate $\norm{\hat{\A} - \A}^2 = O_P(h_n^2)$, and Lemma \ref{lemma:rates} to get $\norm{\X}_F^2 = O_P(p^{1-\delta}n)$. Also we used $\norm{\Delta_j} = 1$ and $\norm{\I_p - \hat{\A}\hat{\A}^T} \leq 2$. For $I_2$, consider
\begin{align}
I_2 &\leq  |n^{-1}s^{-1}\norm{\Delta_j\E}_F^2 - \sigma_j^2| + n^{-1}s^{-1}\norm{\Delta_j\hat{\A}\hat{\A}^T\E}_F^2 \notag\\
&=O_P(n^{-l_{\epsilon}}) + O_P(n^{-1}s^{-1}(\norm{\Delta_j\A\A^T\E}_F^2 + \norm{\Delta_j\A(\hat{\A} - \A)^T\E}_F^2)) \notag\\
&=O_P(n^{-l_{\epsilon}}) + O_P(n^{-1}s^{-1}(\norm{\A^T\E}_F^2 + \norm{(\hat{\A} - \A)^T\E}_F^2)) \notag\\
&= O_P(n^{-l_{\epsilon}}) + O_P(n^{-1}s^{-1}(\norm{\A^T\E}_F^2)), \label{eqn:thm:Cov:I2intermediate}
\end{align}
where we used assumption (D) in arriving at
$|n^{-1}s^{-1}\norm{\Delta_j\E}_F^2 - \sigma_j^2| =
o_P(n^{-l_{\epsilon}})$, and that $\norm{(\hat{\A} - \A)^T\E}_F
\leq \norm{\A^T\E}_F$ for sufficiently large $n$ since
$\norm{\hat{\A} - \A} = o_P(1)$. Consider $\a_i^T\bepsilon_j$
which is the $(i,j)$-th element of $\A^T\E$. We have
\begin{align*}
  E(\a_i^T\bepsilon_j) = 0, \;\;\; \var(\a_i^T\bepsilon_j) = \a_i^T\bSigma_{\bepsilon}\a_i \leq \max_{1\leq j \leq p} \sigma_j^2 = O(1)
\end{align*}
since the $\sigma_j^2$'s are uniformly bounded away from infinity by assumption (M1). Hence each element in $\A^T\E$ is $O_P(1)$, which implies that $$n^{-1}s^{-1}\norm{\A^T\E}_F^2 = O_P(n^{-1}s^{-1} \cdot rn) = O_P(s^{-1}).$$
Hence from (\ref{eqn:thm:Cov:I2intermediate}) we have
\begin{equation}\label{eqn:thm:Cov:I2}
  I_2 = O_P(n^{-l_{\epsilon}} + s^{-1}).
\end{equation}
Assumption (M2) ensures that both $I_1$ and $I_2$ are $o_P(1)$ from (\ref{eqn:thm:Cov:I1}) and (\ref{eqn:thm:Cov:I2}) respectively. From these we can see that
$ I_3 = O_P(I_1^{1/2}) = O_P((p^{1-\delta}s^{-1})^{1/2}h_n)$,
which shows that
\begin{equation}\label{eqn:thm:Cov:Sigmaepsilonhat}
\norm{\hat{\bSigma}_{\bepsilon} - \bSigma_{\bepsilon}} = O_P((p^{1-\delta}s^{-1})^{1/2}h_n).
\end{equation}

Next we consider $\norm{\hat{\bSigma}_{\x} - \bSigma_{\x}} \leq K_1 + K_2 + K_3 + K_4$, where
\begin{align*}
  K_1 &= \norm{\hat{\A}^T\A\tilde{\bSigma}_{\x}\A^T\hat{\A} - \bSigma_{\x}}, \;\;\; K_2 = \norm{\hat{\A}^T\A\tilde{\bSigma}_{\x, \bepsilon}\hat{\A}},\\
  K_3 &= \norm{\hat{\A}^T\tilde{\bSigma}_{\bepsilon, \x}\A^T\hat{\A}}, \;\;\;
  \quad\quad\quad\quad\quad\quad\quad\;\; K_4 = \norm{\hat{\A}^T(\tilde{\bSigma}_{\bepsilon} - \hat{\bSigma}_{\bepsilon})\hat{\A}},
\end{align*}
where $\tilde{\bSigma}_{\x} = n^{-1}\sum_{t=1}^n (\x_t - \bar{\x})(\x_t - \bar{\x})^T$. Now
\begin{align*}
  K_1 &\leq \norm{\hat{\A}^T\A - \I_r} \cdot \norm{\tilde{\bSigma}_{\x}} + \norm{\tilde{\bSigma}_{\x} - \bSigma_{\x}}\\
  &= O_P(\norm{\hat{\A}^T\A - \I_r} \cdot (\norm{\tilde{\bSigma}_{\x} - \bSigma_{\x}} + \norm{\bSigma_{\x}}) + \norm{\tilde{\bSigma}_{\x} - \bSigma_{\x}})\\
  &=O_P(p^{1-\delta}h_n),
\end{align*}
where we used $\norm{\hat{\A}^T\A - \I_r} = \norm{(\hat{\A} - \A)^T\A} = O_P(h_n)$ by Theorem \ref{thm:A}, $\norm{\bSigma_{\x}} = O(p^{1-\delta})$ from assumption in model (\ref{eqn:factormodelWLOG}), and $\norm{\tilde{\bSigma}_{\x} - \bSigma_{\x}} = O_P(p^{1-\delta}n^{-l_x})$ from Lemma \ref{lemma:rates}. Next, using Lemma \ref{lemma:rates} and the fact that $\bSigma_{\x,\bepsilon} = \0$, we have $$K_2 = O_P(K_3) = O_P(\norm{\tilde{\bSigma}_{\x,\bepsilon} - \bSigma_{\x,\bepsilon}}) = O_P(p^{1-\delta/2}n^{-l_{x\epsilon}}).$$
Finally, using Lemma \ref{lemma:rates} again and (\ref{eqn:thm:Cov:Sigmaepsilonhat}),
$$ K_4 = O_P(\norm{\hat{\bSigma}_{\bepsilon} - \bSigma_{\bepsilon}} + \norm{\tilde{\bSigma}_{\bepsilon} - \bSigma_{\bepsilon}}) = O_P((p^{1-\delta}s^{-1})^{1/2}h_n + pn^{-l_{\epsilon}}). $$
Looking at the rates for $K_1$ to $K_4$, and noting assumption (M2) and the definition of $h_n$, we can easily see that
\begin{equation}\label{eqn:thm:Cov:Sigmaxhat}
  \norm{\hat{\bSigma}_{\x} - \bSigma_{\x}} = O_P(p^{1-\delta}h_n).
\end{equation}
From (\ref{eqn:thm:Cov:Sigmayhatnorm}), combining (\ref{eqn:thm:Cov:Sigmaepsilonhat}) and (\ref{eqn:thm:Cov:Sigmaxhat}) and noting assumption (M2), the rate for $\hat{\bSigma}_{\y}$ in the spectral norm is established, and the proof of the theorem completes. $\square$

\vspace{12pt}

{\bf Proof of Theorem \ref{thm:Invcov}.} We first show the rate for $\tilde{\bSigma}_{\y}$. We use the standard inequality
\begin{equation}\label{eqn:thm:invcov:invnorm}
\norm{M_1^{-1} - M_2^{-1}} \leq \frac{\norm{M_2^{-1}}^2 \cdot \norm{M_1 - M_2}}{1 - \norm{(M_1 - M_2)\cdot M_2^{-1}}},
\end{equation}
with $M_1 = \tilde{\bSigma}_{\y}$ and $M_2 = \bSigma_{\y}$. Under assumption (M1) we have $\norm{\bSigma_{\bepsilon}} \asymp 1 \asymp \norm{\bSigma_{\bepsilon}^{-1}}$, so that
\begin{align*}
\norm{\bSigma_{\y}^{-1}} &= \norm{\bSigma_{\bepsilon}^{-1} - \bSigma_{\bepsilon}^{-1}\A(\bSigma_{\x}^{-1} + \A^T\bSigma_{\bepsilon}^{-1}\A)^{-1}\A^T\bSigma_{\bepsilon}^{-1}}\\
&= O(\norm{\bSigma_{\bepsilon}^{-1}} + \norm{\bSigma_{\bepsilon}^{-1}}^2 \cdot \norm{(\bSigma_{\x}^{-1} + \A^T\bSigma_{\bepsilon}^{-1}\A)^{-1}})\\
&=O(1),
\end{align*}
where we also used
\begin{align}
\norm{(\bSigma_{\x}^{-1} + \A^T\bSigma_{\bepsilon}^{-1}\A)^{-1}}
&\leq \norm{(\A^T\bSigma_{\bepsilon}^{-1}\A)^{-1}} \notag\\ &=
\lambda_{\max}\{({\A}^T{\bSigma}_{\bepsilon}^{-1}{\A})^{-1}\} =
\lambda_{\min}^{-1}({\A}^T{\bSigma}_{\bepsilon}^{-1}{\A})  = O(1),
\label{eqn:thm:invcov:inq1}
\end{align}
since the eigenvalues of $\bSigma_{\bepsilon}^{-1}$ are of constant
order by assumption (M1), and $\norm{\A}_{\min} = 1$ with $\A\x \neq
\0$ for any $\x$ since $\A$ is of full rank with $p > r$, so that
\begin{align}\label{eqn:thm:invcov:norminequality}
\lambda_{\min}(\A^T\bSigma_{\bepsilon}^{-1}\A) = \min_{\x \neq \0}
\frac{\x^T\A^T\bSigma_{\bepsilon}^{-1}\A\x}{\norm{\x}} &\geq
\min_{\y \neq \0} \frac{\y^T\bSigma_{\bepsilon}^{-1}\y}{\norm{\y}}
\cdot \min_{\x
\neq \0}\frac{\norm{\A\x}}{\norm{\x}} \notag\\
&=\lambda_{\min}(\bSigma_{\bepsilon}^{-1})\cdot\norm{\A}_{\min}
\asymp 1.
\end{align}
Then by (\ref{eqn:thm:invcov:invnorm}) together with Theorem
\ref{thm:Cov} that $\norm{\tilde{\bSigma}_{\y} - \bSigma_{\y}} =
O_P(p^{1-\delta}h_n) = o_P(1)$, we have
$$ \norm{\tilde{\bSigma}_{\y}^{-1} - \bSigma_{\y}^{-1}} = \frac{O(1)^2 \cdot O_P(p^{1-\delta}h_n)}{1-o_P(1)} = O_P(p^{1-\delta}h_n), $$
which is what we need to show.

 Now we show the rate for $\hat{\bSigma}_{\y}^{-1}$. Using the Sherman-Morrison-Woodbury formula, we have $\norm{\hat{\bSigma}_{\y}^{-1} - \bSigma_{\y}^{-1}} \leq \sum_{j=1}^6 K_j$, where
\begin{equation}\label{eqn:thm:invcov:sigmayhatinv}
\begin{split}
    K_1 &= \norm{\hat{\bSigma}_{\bepsilon}^{-1} - \bSigma_{\bepsilon}^{-1}}, \\
    K_2 &= \norm{(\hat{\bSigma}_{\bepsilon}^{-1} - \bSigma_{\bepsilon}^{-1})\hat{\A}(\hat{\bSigma}_{\x}^{-1} + \hat{\A}^T\hat{\bSigma}_{\bepsilon}^{-1}\hat{\A})^{-1}\hat{\A}^T\hat{\bSigma}_{\bepsilon}^{-1}}, \\
    K_3 &= \norm{\bSigma_{\bepsilon}^{-1}\hat{\A}(\hat{\bSigma}_{\x}^{-1} + \hat{\A}^T\hat{\bSigma}_{\bepsilon}^{-1}\hat{\A})^{-1}\hat{\A}^T(\hat{\bSigma}_{\bepsilon}^{-1} - \bSigma_{\bepsilon}^{-1})}, \\
    K_4 &= \norm{\bSigma_{\bepsilon}^{-1}(\hat{\A} - \A)(\hat{\bSigma}_{\x}^{-1} + \hat{\A}^T\hat{\bSigma}_{\bepsilon}^{-1}\hat{\A})^{-1}\hat{\A}^T \bSigma_{\bepsilon}^{-1}}, \\
    K_5 &= \norm{\bSigma_{\bepsilon}^{-1}\A(\hat{\bSigma}_{\x}^{-1} + \hat{\A}^T\hat{\bSigma}_{\bepsilon}^{-1}\hat{\A})^{-1}(\hat{\A} - \A)^T \bSigma_{\bepsilon}^{-1}}, \\
    K_6 &= \norm{\bSigma_{\bepsilon}^{-1}\A\{ (\hat{\bSigma}_{\x}^{-1} + \hat{\A}^T\hat{\bSigma}_{\bepsilon}^{-1}\hat{\A})^{-1} - (\bSigma_{\x}^{-1} + \A^T\bSigma_{\bepsilon}^{-1}\A)^{-1} \}\A^T\bSigma_{\bepsilon}^{-1}}.
\end{split}
\end{equation}
First, we have $\norm{\bSigma_{\bepsilon}^{-1}} = O(1)$ as before by assumption (M1). Next,
$$ K_1 = \max_{1\leq j \leq k}|\hat{\sigma}_j^{-2} - \sigma_j^{-2}| = \frac{\max_{1 \leq j \leq k}|\hat{\sigma}_j^2 - \sigma_j^2|}{\min_{1 \leq j \leq k} \sigma_j^2(\sigma_j^2 + O_P((p^{1-\delta}s^{-1})^{1/2} h_n ))} = O_P((p^{1-\delta}s^{-1})^{1/2} h_n ), $$
where we used (\ref{eqn:thm:Cov:Sigmaepsilonhat}) and assumptions
(M1) and (M2). From these, we have
\begin{equation}\label{eqn:thm:invcov:Sigmaepsilonhat}
\norm{\hat{\bSigma}_{\bepsilon}^{-1}} = O_P(1).
\end{equation}
Also, like (\ref{eqn:thm:invcov:inq1}),
\begin{align}
  \norm{(\hat{\bSigma}_{\x} + \hat{\A}^T\hat{\bSigma}_{\bepsilon}^{-1}\hat{\A})^{-1}} = O_P(1), \label{eqn:thm:invcov:inq2}
\end{align}
noting (\ref{eqn:thm:Cov:Sigmaepsilonhat}) and assumption (M1). With
these rates and noting that $\norm{\A} = \norm{\hat{\A}} = 1$ and
$\norm{\hat{\A} - \A} = O_P(h_n)$ from Theorem \ref{thm:A}, we
have from (\ref{eqn:thm:invcov:sigmayhatinv}) that
\begin{equation}\label{eqn:thm:invcov:sigmayhatinv1}
  \begin{split}
    \norm{\hat{\bSigma}_{\y}^{-1} - \bSigma_{\y}^{-1}} &= O_P((p^{1-\delta}s^{-1})^{1/2}h_n + h_n)\\
    &\;\;\; + O_P(\norm{(\hat{\bSigma}_{\x}^{-1} + \hat{\A}^T\hat{\bSigma}_{\bepsilon}^{-1}\hat{\A})^{-1} - (\bSigma_{\x}^{-1} + \A^T\bSigma_{\bepsilon}^{-1}\A)^{-1}}),
  \end{split}
\end{equation}
where the last term is contributed from $K_6$. Using (\ref{eqn:thm:invcov:inq1}) and (\ref{eqn:thm:invcov:inq2}), and the inequality
$\norm{M_1^{-1} - M_2^{-1}} = O_P(\norm{M_1^{-1}}\cdot \norm{M_1-M_2}\cdot \norm{M_2^{-1}})$, the rate for this term can be shown to be $O_P(L_1 + L_2)$, where
\begin{align*}
  L_1 = \norm{\hat{\bSigma}_{\x}^{-1} - \bSigma_{\x}^{-1} }, \;\;\; L_2 = \norm{\hat{\A}^T\hat{\bSigma}_{\bepsilon}^{-1}\hat{\A} - \A^T\bSigma_{\bepsilon}^{-1}\A}.
\end{align*}
Consider $\norm{\bSigma_{\x}^{-1}} \leq \norm{\bSigma_{\x}^{-1}} = \lambda_{\min}^{-1}(\bSigma_{\x}) = O(p^{-(1-\delta)})$ by assumption in model (\ref{eqn:factormodelWLOG}). With this and (\ref{eqn:thm:Cov:Sigmaxhat}), substituting $M_1 = \hat{\bSigma}_{\x}$ and $M_2 = \bSigma_{\x}$ into (\ref{eqn:thm:invcov:invnorm}), we have
\begin{equation}\label{eqn:thm:invcov:L1}
L_1 = \frac{O(p^{-(2-2\delta)}) \cdot O_P(p^{1-\delta}h_n)}{1-O_P(p^{1-\delta}h_n \cdot p^{-(1-\delta)})} = \frac{O_P(p^{-(1-\delta)}h_n)}{1-o_P(1)} = O_P(p^{-(1-\delta)}h_n).
\end{equation}
For $L_2$, using $\norm{\hat{\A}} = \norm{\A} = 1$, $\norm{\hat{\A} - \A} = O_P(h_n)$ from Theorem \ref{thm:A}, the rate for $K_1$ shown before and (\ref{eqn:thm:invcov:Sigmaepsilonhat}), we have
\begin{align}\label{eqn:thm:invcov:L2}
  L_2 = O_P(\norm{\hat{\A} - \A} \cdot \norm{\hat{\bSigma}_{\bepsilon}^{-1}} + \norm{\hat{\bSigma}_{\bepsilon}^{-1} - \bSigma_{\bepsilon}^{-1}}) = O_P(h_n + (p^{1-\delta}s^{-1})^{1/2}h_n).
\end{align}
Hence, from (\ref{eqn:thm:invcov:sigmayhatinv1}), together with (\ref{eqn:thm:invcov:L1}) and (\ref{eqn:thm:invcov:L2}), we have $$\norm{\hat{\bSigma}_{\y}^{-1} - \bSigma_{\y}^{-1}} = O_P((1+(p^{1-\delta}s^{-1})^{1/2})h_n),$$ which completes the proof of the theorem. $\square$

\vspace{12pt}

{\bf Proof of Theorem \ref{thm:stationaryfactors}.} The idea of the proof is similar to that for Theorem \ref{thm:A} for the simple procedure. We want to find the order of the eigenvalues of the matrix $\L$ first.

 From model (\ref{eqn:statfactormodelWLOG}), we have for $i=1,2$,
 \begin{equation}\label{eqn:thm:Ans1:inq1}
 \begin{split}
 \norm{\bSigma_{ii}(k)}^2 &\asymp p^{2-2\delta_i} \asymp \norm{\bSigma_{ii}(k)}_{\min}^2,\\
 \norm{\bSigma_{12}(k)}^2 &= O(p^{2-\delta_1 - \delta_2}) = \norm{\bSigma_{21}(k)}^2,\\
 \norm{\bSigma_{i\bepsilon}(k)}^2 &= O(p^{2-\delta_i}).
 \end{split}
 \end{equation}
  We want to find the lower bounds of the order of the $r_1$-th largest eigenvalue, as well as the smallest non-zero eigenvalue of $\L$. We
first note that
\begin{align}
\bSigma_{\y}(k) &= \A_1\bSigma_{11}(k)\A_1^T + \A_1\bSigma_{12}(k)\A_2^T + \A_2\bSigma_{21}(k)\A_1^T \notag\\
&\;\;\;+\A_2\bSigma_{22}(k)\A_2^T  + \A_1\bSigma_{1\bepsilon}(k) + \A_2\bSigma_{2\bepsilon}(k), \label{eqn:thm:Ans1:Sy}
\end{align}
and hence
\begin{align}\label{eqn:thm:Ans1:Lpoundexpansion}
  \L = \A_1\W_1\A_1^T + \A_2\W_2\A_2^T + \text{cross terms},
\end{align}
where $\W_1$ (with size $r_1 \times r_1$) and $\W_2$ (with size $r_2 \times r_2$) are positive semi-definite matrices defined
by
\begin{align*}
  \W_1 &= \sum_{k=1}^{k_0} \Big\{ \bSigma_{11}(k)\bSigma_{11}(k)^T + \bSigma_{12}(k)\bSigma_{12}(k)^T + \bSigma_{1\bepsilon}(k)\bSigma_{1\bepsilon}(k)^T\Big\},\\
  \W_2 &= \sum_{k=1}^{k_0} \Big\{ \bSigma_{21}(k)\bSigma_{21}(k)^T + \bSigma_{2\bepsilon}(k)\bSigma_{2\bepsilon}(k)^T + \bSigma_{22}(k)\bSigma_{22}(k)^T
  \Big\}.
\end{align*}
From $\W_1$, by (\ref{eqn:thm:Ans1:inq1}) and that $\norm{\bSigma_{1\bepsilon}(k)} = o(p^{1-\delta_1})$, we have the order of the $r_1$ eigenvalues for $\W_1$ is all $p^{2-2\delta_1}$. Then the $r_1$-th largest eigenvalue of $\L$ is of order $p^{2-2\delta_1}$ since the term $\bSigma_{11}(k)\bSigma_{11}(k)^T$ has the largest order at $p^{2-2\delta_1}$. We
 write
\begin{equation}\label{eqn:Lr1eigenvalue}
p^{2-2\delta_1} = O(\lambda_{r_1}(\L)),
\end{equation}
where $\lambda_i(M)$ represents the $i$-th largest eigenvalue of the square matrix $M$.

For the smallest non-zero eigenvalue of $\W_2$, since $\norm{\bSigma_{2\bepsilon}(k)} = o(p^{1-\delta_2})$, it
is contributed either from the term $\bSigma_{22}(k)\bSigma_{22}(k)^T$ or $\bSigma_{21}(k)\bSigma_{21}(k)^T$ in $\W_2$, and has order $p^{2-2\delta_2}$ if $\norm{\bSigma_{21}(k)}_{\min} = O(p^{2-2\delta_2})$, and $p^{2-c}$ in general if $\norm{\bSigma_{21}(k)}_{\min} \asymp p^{2-c}$, with $\delta_1 + \delta_2 \leq c \leq 2\delta_2$. Hence
\begin{equation}\label{eqn:Lr1r2eigenvalue}
\begin{split}
p^{2-2\delta_2} &= O(\lambda_{r_1+r_2}(\L)), \;\; \text{ if } \norm{\bSigma_{21}(k)}_{\min} = o(p^{2-2\delta_2}),\\
p^{2-c} &= O(\lambda_{r_1+r_2}(\L)), \;\; \text{ if } \norm{\bSigma_{21}(k)}_{\min} \asymp p^{2-c}, \; \delta_1 + \delta_2 \leq c < 2\delta_2.
\end{split}
\end{equation}

Now we can write $\L = (\A \; \B)\D(\A \; \B)^T$, where $\B$ is the orthogonal complement of $\A$, and
with $\D_{1}$ containing the $r_1$ largest eigenvalues of $\L$ and $\D_{2}$ the next $r_2$ largest,
\begin{equation}\label{eqn:thm:Ans1:Dx}
  \D = \left(
           \begin{array}{ccc}
             \D_{1} & \0 & \0 \\
           \0 & \D_{2} & \0 \\
           \0 & \0 & \0
           \end{array}
         \right),
\end{equation}
with $p^{2-2\delta_1} = O(\lambda_{\min}(\D_{1}))$ by (\ref{eqn:Lr1eigenvalue}), and $p^{2-c} = O(\lambda_{\min}(\D_{2}))$ by (\ref{eqn:Lr1r2eigenvalue}), $\delta_1 + \delta_2 \leq c \leq 2\delta_2$.

Similar to the proof of Theorem \ref{thm:A}, we define $\E_{\L} =\tilde{\L} - \L$. Then (\ref{eqn:thm:A:E_Lbound}) holds, and
\begin{equation}\label{eqn:thm:Ans1:Synorm}
  \norm{\bSigma_{\y}(k)} = O(p^{1-\delta_1}),
\end{equation}
using (\ref{eqn:thm:Ans1:Sy}) and $\norm{\bSigma_{11}(k)} = O(p^{1-\delta_1})$ from (\ref{eqn:thm:Ans1:inq1}).
Also,
\begin{align}
\norm{\tilde{\bSigma}_{\y}(k) - \bSigma_{\y}(k)} &= O_P(\norm{\tilde{\bSigma}_{11}(k) - \bSigma_{11}(k)} + \norm{\tilde{\bSigma}_{12}(k) - \bSigma_{12}(k)} \notag\\
&\;\;\; + \norm{\tilde{\bSigma}_{22}(k) - \bSigma_{22}(k)} + \norm{\tilde{\bSigma}_{1\bepsilon}(k) - \bSigma_{1\bepsilon}(k)} \notag\\
&\;\;\; + \norm{\tilde{\bSigma}_{2\bepsilon}(k) - \bSigma_{2\bepsilon}(k)} + \norm{\tilde{\bSigma}_{\bepsilon}(k) - \bSigma_{\bepsilon}(k)} ) \notag\\
&=O_P(p^{1-\delta_1}n^{-l_1} + p^{1-\delta_2}n^{-l_{2}} + p^{1-\delta_1/2-\delta_2/2}n^{-l_{12}} \notag\\
&\;\;\; + p^{1-\delta_1/2}n^{-l_{1\epsilon}} + p^{1-\delta_2/2}n^{-l_{2\epsilon}} + pn^{-l_{\epsilon}}) = O_P(p^{1-\delta_1}\omega_1), \label{eqn:thm:Ans1:Sytildenorm}
\end{align}
where we used condition (D') in section \ref{subsect:theories} to derive the following rates like those in Lemma \ref{lemma:rates} (proofs thus omitted):
\begin{equation}\label{eqn:thm:Ans1:rates}
\begin{split}
\norm{\tilde{\bSigma}_{11}(k) - \bSigma_{11}(k)} &= O_P(p^{1-\delta_1}n^{-l_{1}}), \;\;\;\;\; \norm{\tilde{\bSigma}_{12}(k) - \bSigma_{12}(k)} = O_P(p^{1-\delta_1/2 - \delta_2/2}n^{- l_{12}}),\\
\norm{\tilde{\bSigma}_{22}(k) - \bSigma_{22}(k)} &= O_P(p^{1-\delta_2}n^{-l_{2}}), \;\;\;\;\;
\norm{\tilde{\bSigma}_{1\bepsilon}(k) - \bSigma_{1\bepsilon}(k)} = O_P(p^{1-\delta_1/2}n^{-l_{1\epsilon}}),\\
\norm{\tilde{\bSigma}_{2\bepsilon}(k) - \bSigma_{2\bepsilon}(k)} &= O_P(p^{1-\delta_2/2}n^{-l_{2\epsilon}}), \;\;\;\,
\norm{\tilde{\bSigma}_{\bepsilon}(k) - \bSigma_{\bepsilon}(k)} = O_P(pn^{-l_{\epsilon}}).
\end{split}
\end{equation}

We form $\hat{\A}_1$ with the first $r_1$ unit eigenvectors corresponding to the $r_1$ largest eigenvalues, i.e. the eigenvalues in $\D_{1}$. Now, we have
\begin{align*}
 \norm{\E_{\L}} &= O_P(\norm{\tilde{\bSigma}_{\y}(k) - \bSigma_{\y}(k)} \cdot (\norm{\bSigma_{\y}(k)} + \norm{\tilde{\bSigma}_{\y}(k) - \bSigma_{\y}(k)})) = O_P(p^{2-2\delta_1}\omega_1) \\
 &= o_P(p^{2-2\delta_1}) = O_P(\lambda_{\min}(\D_{1})), \text{ hence for } n \text{ large enough, }\\
 \norm{\E_{\L}} &\leq  \frac{1}{5}\sep\bigg(\D_{1}, \left(
                                      \begin{array}{cc}
                                        \D_{2} & \0 \\
                                        \0 & \0 \\
                                      \end{array}
                                    \right)
\bigg),
\end{align*}
where the second equality is from (\ref{eqn:thm:Ans1:Synorm}) and (\ref{eqn:thm:Ans1:Sytildenorm}), the third is from noting that $\omega_1 = o(1)$, and the last is from (\ref{eqn:Lr1eigenvalue}). Hence, we can use Lemma \ref{lemma:matrixperturbation} and arguments similar to the proof of Theorem \ref{thm:A} to conclude that
\begin{align*}
\norm{\hat{\A}_1 - \A_1} = O_P\bigg(\norm{\E_{\L}}/{\sep\bigg(\D_{1}, \left(
                                      \begin{array}{cc}
                                        \D_{2} & \0 \\
                                        \0 & \0 \\
                                      \end{array}
                                    \right)
\bigg)}\bigg) = O_P(\omega_1).
\end{align*}
Similarly, depending on the order of $\norm{\bSigma_{21}(k)}_{\min}$, we have
\begin{align*}
 \norm{\E_{\L}} &= O_P(p^{2-2\delta_1}\omega_1)\leq  \frac{1}{5}\sep\bigg(\D_{2}, \left(
                                      \begin{array}{cc}
                                        \D_{1} & \0 \\
                                        \0 & \0 \\
                                      \end{array}
                                    \right)
\bigg)
\asymp p^{2-c},
\end{align*}
since we assumed $p^{c-2\delta_1}\omega_1 = o(1)$ for $\delta_1 + \delta_2 \leq c \leq 2\delta_2$. Hence
$$ \norm{\hat{\A}_2 - \A_{2}} = O_P\bigg(\norm{\E_{\L}}/{\sep\bigg(\D_{2}, \left(
                                      \begin{array}{cc}
                                        \D_{1} & \0 \\
                                        \0 & \0 \\
                                      \end{array}
                                    \right)
\bigg)}\bigg) = O_P(p^{c - 2\delta_1}\omega_1) = O_P(\omega_2). $$
This completes the proof for the simple procedure.

\vspace{12pt}

For the two-step procedure,
denote $\y_t^* = (\I_p - \hat{\A}_1\hat{\A}_1^T)\y_t$, and define $\E_{\L^*} = \tilde{\L}^* - \L^*$. Note that
\begin{equation}\label{eqn:thm:Ans2:EL*}
 \norm{\E_{\L^*}} = \norm{\tilde{\L}^* - \L^*} \leq \sum_{k=1}^{k_0} \Big\{ \norm{\tilde{\bSigma}_{\y^*}(k) - \bSigma_{\y^*}(k)}^2 + \norm{\bSigma_{\y^*}(k)} \cdot \norm{\tilde{\bSigma}_{\y^*}(k) - \bSigma_{\y^*}(k)} \Big\},
\end{equation}
where
\begin{align*}
  \tilde{\bSigma}_{\y^*}(k) &= (\I_p - \hat{\A}_1\hat{\A}_1^T) \tilde{\bSigma}_{\y}(k) (\I_p - \hat{\A}_1\hat{\A}_1^T),\\
  \bSigma_{\y^*}(k) &= (\I_p - \A_1\A_1^T)\bSigma_{\y}(k)(\I_p - \A_1\A_1^T) = \A_2\bSigma_{22}(k)\A_2^T
  + \A_2\bSigma_{2\bepsilon}(k)(\I_p - \A_1\A_1^T),\\
  \tilde{\L}^* &= \sum_{k=1}^{k_0} \tilde{\bSigma}_{\y^*}(k)\tilde{\bSigma}_{\y^*}(k)^T, \;\;\; \L^*= \sum_{k=1}^{k_0} \bSigma_{\y^*}(k)\bSigma_{\y^*}(k)^T,
\end{align*}
with $\hat{\A}_1$ being the estimator from the simple procedure, so that $\norm{\hat{\A}_1 - \A_1} = O_P(\omega_1)$ from previous result. We write
$$\L^* = \A_2\Q_2\D_2^*\Q_2^T\A_2^T, \text{ so that } \L^*\A_2\Q_2 = \A_2\Q_2\D_2^*,$$
and like section \ref{subsect:estimation}, we take $\A_2\Q_2$ as the $\A_2$ to be used in our inference.

The idea of the proof is to find the rates of $\norm{\E_{\L^*}}$ and the eigenvalues in $\D_{2}^*$ and use the arguments similar to the proof for the simple procedure to get the rate for $\norm{\check{\A}_2 - \A_2}$.

First, with the assumption that $\norm{\bSigma_{22}(k)} \asymp p^{1-\delta_2} \asymp \norm{\bSigma_{22}(k)}_{\min}$ and $\norm{\bSigma_{2\bepsilon}(k)} = o(p^{1-\delta_2})$, all the eigenvalues in $\D_2^*$ have order $p^{2-2\delta_2}$.

We need to find $\norm{\E_{\L^*}}$. It is easy to show that
\begin{equation}\label{eqn:thm:Ans2:A1A2rates}
\begin{split}
  \norm{\hat{\A}_1\hat{\A}_1^T\A_2} &= \norm{\hat{\A}_1(\hat{\A}_1 - \A_1)^T\A_2} \leq \norm{\hat{\A}_1 - \A_1} = O_P(\omega_1),\\
  \norm{(\I_p - \hat{\A}_1\hat{\A}_1^T)\A_1} &= \norm{(\A_1 - \hat{\A}_1) - \hat{\A}_1(\hat{\A}_1 - \A_1)^T\A_1 } = O_P(\omega_1).
\end{split}
\end{equation}

Writing $\hat{\H}_1 = \I_p - \hat{\A}_1\hat{\A}_1^T$, we can decompose $\tilde{\bSigma}_{\y^*}(k) - \bSigma_{\y^*}(k) = \sum_{l=1}^9 I_l$, where
\begin{align*}
  I_1 &= \hat{\H}_1\A_1\tilde{\bSigma}_{11}(k)\A_1^T\hat{\H}_1,\;\;   I_2 = \hat{\H}_1\A_1\tilde{\bSigma}_{12}(k)\A_2^T\hat{\H}_1,\\ I_3 &= \hat{\H}_1\A_1\tilde{\bSigma}_{1\bepsilon}(k)\hat{\H}_1, \quad\quad I_4 = \hat{\H}_1\A_2\tilde{\bSigma}_{21}(k)\A_1^T\hat{\H}_1,\\
  I_5 &= \hat{\H}_1\A_2\tilde{\bSigma}_{22}(k)\A_2^T\hat{\H}_1 - \A_2\bSigma_{22}(k)\A_2^T, \;\;\; I_6 = \hat{\H}_1\A_2\tilde{\bSigma}_{2\bepsilon}(k)\hat{\H}_1 - \A_2\bSigma_{2\bepsilon}(k)\H_1,\\
  I_7 &= \hat{\H}_1\tilde{\bSigma}_{\bepsilon 1}(k)\A_1^T\hat{\H}_1, \;\;\, I_8 = \hat{\H}_1\tilde{\bSigma}_{\bepsilon 2}(k)\A_2^T\hat{\H}_1, \;\;\;\, I_9 = \hat{\H}_1\tilde{\bSigma}_{\bepsilon}(k)\hat{\H}_1.
\end{align*}
Using (\ref{eqn:thm:Ans1:rates}), (\ref{eqn:thm:Ans2:A1A2rates}) and the assumptions in model (\ref{eqn:statfactormodelWLOG}), we can see that
\begin{align*}
  \norm{I_1} &= O_P(\norm{\hat{\H}_1\A_1}^2\cdot\norm{\tilde{\bSigma}_{11}(k)}) = O_P(p^{1-\delta_1}\omega_1^2),\\
 \norm{I_2} &= O_P(\norm{\hat{\H}_1\A_1} \cdot \norm{\tilde{\bSigma}_{12}(k)}) = O_P(p^{1-\delta_1/2-\delta_2/2}\omega_1) = \norm{I_4},\\
 \norm{I_7} &= O_P(\norm{I_3}) = O_P(\norm{\hat{\H}_1\A_1} \cdot (\norm{\tilde{\bSigma}_{1\bepsilon}(k) - \bSigma_{1\bepsilon}(k)} + \norm{\bSigma_{1\bepsilon}(k)})) = O_P(p^{1-\delta_1}\omega_1),\\
 \norm{I_6} &= O_P(\norm{\tilde{\bSigma}_{2\bepsilon}(k) - \bSigma_{2\bepsilon}(k)} + \norm{\bSigma_{2\bepsilon}(k)}\cdot\norm{\hat{\A}_1 - \A_1}) = O_P(p^{1-\delta_2/2}(n^{-l_{2\epsilon}} + p^{-\delta_2/2}\omega_1)),\\ \norm{I_8} &= O_P(\norm{\tilde{\bSigma}_{\bepsilon 2}(k)}) = O_P(p^{1-\delta_2/2}n^{-l_{2\epsilon}}), \;\;\; \norm{I_9} = O_P(\norm{\tilde{\bSigma}_{\bepsilon}}) = O_P(pn^{-l_{\epsilon}}),\\
 \norm{I_5} &= O_P(\norm{\hat{\A}_1\hat{\A}_1^T\A_2}\cdot\norm{\tilde{\bSigma}_{22}(k)} + \norm{\tilde{\bSigma}_{22}(k) - \bSigma_{22}(k)}) = O_P(p^{1-\delta_2}(\omega_1 + n^{-l_2})).
\end{align*}
Hence, we have
\begin{align}
  \norm{\tilde{\bSigma}_{\y^*}(k) - \bSigma_{\y^*}(k)} &= O_P(p^{1-\delta_1}\omega_1^2 + p^{1-\delta_1/2-\delta_2/2}\omega_1 + p^{1-\delta_1}\omega_1 \notag\\
&\;\;\; + p^{1-\delta_2/2}(n^{-l_{2\epsilon}} + p^{-\delta_2/2}\omega_1) + pn^{-l_{\epsilon}} + p^{1-\delta_2}(\omega_1 + n^{-l_2})) \notag\\
&= O_P(p^{1-\delta_1}\omega_1). \label{eqn:thm:Ans2:Sigmatildey*}
\end{align}
We also have
\begin{equation}\label{eqn:thm:Ans2:Sigmay*norm}
  \norm{\bSigma_{\y^*}(k)} = O_P(\norm{\bSigma_{22}(k)} + \norm{\bSigma_{2\bepsilon}(k)}) = O_P(p^{1-\delta_2}).
\end{equation}

Hence, with (\ref{eqn:thm:Ans2:Sigmatildey*}) and (\ref{eqn:thm:Ans2:Sigmay*norm}), (\ref{eqn:thm:Ans2:EL*}) becomes
\begin{equation}\label{eqn:thm:Ans2:EL*norm}
  \norm{\E_{\L^*}} = O_P(p^{2-\delta_1 - \delta_2}\omega_1).
\end{equation}

With the order of eigenvalues in $\D_2^*$ being $p^{2-2\delta_2}$ and noting (\ref{eqn:thm:Ans2:Sigmatildey*}), we can use Lemma \ref{lemma:matrixperturbation} and the arguments similar to those in the proof of Theorem \ref{thm:A} to get
$$ \norm{\check{\A}_2 - \A_2} = O_P(\norm{\E_{\L^*}}/\sep(\D_2^*, \0)) = O_P(p^{\delta_2-\delta_1}\omega_1), $$
and the proof of the theorem completes. $\square$

\vspace{12pt}

{\bf Proof of Theorem \ref{thm:comparefactormagnitude}.} We have $\hat{\x}_t = \hat{\A}^T\y_t = \hat{\A}^T\A\x_t + (\hat{\A} - \A)^T\bepsilon_t + \A^T\bepsilon_t$. With $\A = (\A_1 \; \A_2)$ and $\x_t = (\x_{1t}^T \; \x_{2t}^T)^T$, we have
\begin{align*}
  \hat{\x}_{1t} &= \hat{\A}_1^T\A_1\x_{1t} + \hat{\A}_1^T\A_2\x_{2t} + (\hat{\A}_1 - \A_1)^T\bepsilon_t + \A_1^T\bepsilon_t,\\
  \hat{\x}_{2t} &= \hat{\A}_2^T\A_1\x_{1t} + \hat{\A}_2^T\A_2\x_{2t} + (\hat{\A}_2 - \A_2)^T\bepsilon_t + \A_2^T\bepsilon_t.
\end{align*}
We first note that for $i=1,2$, $\norm{(\hat{\A}_i - \A_i)^T\bepsilon_t} = O_P(\A_i^T\bepsilon_t) = O_P(1)$ since $\norm{\hat{\A}_i - \A_i} = o_P(1)$ and $\A_i^T\bepsilon_t$ are $r_i$ $O_P(1)$ random variables. Then
\begin{align*}
  \norm{\hat{\x}_{1t}} = \norm{\hat{\x}_{1t}}_F &\geq \norm{\hat{\A}_1^T\A_1\x_{1t}}_F - \norm{\hat{\A}_1^T\A_2\x_{2t}}_F + O_P(1)\\
                       &\geq \norm{\hat{\A}_1^T\A_1}_{\min} \cdot \norm{\x_{1t}}_F - \norm{\hat{\A}_1^T\A_2}\cdot \norm{\x_{2t}} + O_P(1)\\
&\geq \norm{\hat{\A}_1}_{\min} \cdot \norm{\A_1}_{\min} \cdot \norm{\x_{1t}} - o_P(\norm{\x_{2t}}) + O_P(1)\\
&\asymp_P \norm{\x_{1t}} \asymp p^{\frac{1-\delta_1}{2}},
\end{align*}
where $\norm{M}_F$ denotes the Frobenius norm of the matrix $M$, and we used the inequality $\norm{\A\B}_F \geq \norm{\A}_{\min} \cdot \norm{\B}_F$. Finally, with similar arguments,
\begin{align*}
  \norm{\hat{\x}_{2t}} &\leq \norm{\x_{2t}} + \norm{\hat{\A}_2^T\A_2} \cdot \norm{\x_{1t}} + O_P(1)\\
&= O_P(p^{\frac{1-\delta_2}{2}}) + O_P(\norm{\hat{\A}_2 - \A_2}\cdot p^{\frac{1-\delta_1}{2}})\\
&= O_P(p^{\frac{1-\delta_2}{2}}) + o_P(p^{\frac{1-\delta_1}{2}})\\
&= o_P(p^{\frac{1-\delta_1}{2}}),
\end{align*}
which establishes the claim of the theorem. $\square$

\vspace{12pt}

{\bf Proof of Theorem \ref{thm:statcov}.} We can easily use the decomposition in (\ref{eqn:thm:A:Sigmaytildenorm}) again for model (\ref{eqn:statfactormodelWLOG}) to arrive at
\begin{align*}
\norm{\tilde{\bSigma}_{\y} - \bSigma_{\y}} &= O_P(p^{1-\delta_1}n^{-l_1} + p^{1-\delta_2}n^{-l_2} + p^{1-\delta_1/2}n^{-l_{1\epsilon}} + p^{1-\delta_2/2}n^{-l_{2\epsilon}} + pn^{-l_{\epsilon}}),\\
&=O_P(p^{1-\delta_1}\omega_1),
\end{align*}
where we used assumption (D'), and arguments like those in Lemma \ref{lemma:rates} to arrive at $\norm{\tilde{\bSigma}_{\x} - \bSigma_{\x}} = O_P(p^{1-\delta_1}n^{-l_1} + p^{1-\delta_2}n^{-l_2})$, $\norm{\tilde{\bSigma}_{\x, \bepsilon} - \bSigma_{\x, \bepsilon}} = O_P(p^{1-\delta_1/2}n^{-l_{1\epsilon}} + p^{1-\delta_2/2}n^{-l_{2\epsilon}})$ and $\norm{\tilde{\bSigma}_{\bepsilon} - \bSigma_{\bepsilon}} = O_P(pn^{-l_{\epsilon}})$.

Now consider $\norm{\hat{\bSigma}_{\y} - \bSigma_{\y}}$. The proof for $\norm{\check{\bSigma}_{\y} - \bSigma_{\y}}$ follows exactly the same lines by replacing $\hat\A$ with $\check\A$ and is thus omitted. It can be decomposed like that in (\ref{eqn:thm:Cov:Sigmayhatnorm}). Hence we need to consider $\norm{\hat{\bSigma}_{\bepsilon} - \bSigma_{\bepsilon}} = \max_{1 \leq j \leq p} |\hat{\sigma}_j^2 - \sigma_j^2| \leq I_1 + I_2 + I_3$, where $$I_1 = O_P(p^{1-\delta_1}s^{-1}\norm{\hat{\A} - \A}^2),$$
which used decomposition in (\ref{eqn:thm:Cov:I1}), and $\norm{\X}_F^2 = O_P(p^{1-\delta_1}nr_1 + p^{1-\delta_2}nr_2) = O_P(p^{1-\delta_1}n)$;
$$ I_2 = O_P(n^{-l_{\epsilon}} + s^{-1}), $$
where derivation is similar to that in (\ref{eqn:thm:Cov:I2intermediate}) and thereafter. Also, $I_3 = O_P(I_1^{1/2}) = O_P((p^{1-\delta_1}s^{-1})^{1/2}\norm{\hat{\A} - \A})$. Thus, with assumption (M2)', we see that
\begin{equation}\label{eqn:thm:statcov:Sigmahatepsilonnorm}
  \norm{\hat{\bSigma}_{\bepsilon} - \bSigma_{\bepsilon}} = O_P((p^{1-\delta_1}s^{-1})^{1/2}\norm{\hat{\A} - \A}).
\end{equation}

For $\norm{\hat{\bSigma}_{\x} - \bSigma_{\x}}$, we use the decomposition like that in the proof of Theorem \ref{thm:Cov}, and noting that $\norm{\bSigma_{\x}} = O(p^{1-\delta_1} + p^{1-\delta_2}) = O(p^{1-\delta_1})$, to arrive at
\begin{equation}\label{eqn:thm:statcov:Sigmahatxnorm}
  \norm{\hat{\bSigma}_{\x} - \bSigma_{\x}} = O_P(p^{1-\delta_1}\norm{\hat{\A} - \A}).
\end{equation}
Hence noting assumption (M2)' again and combining (\ref{eqn:thm:statcov:Sigmahatepsilonnorm}) and (\ref{eqn:thm:statcov:Sigmahatxnorm}), we see that
$$ \norm{\hat{\bSigma}_{\y} - \bSigma_{\y}} = O_P(p^{1-\delta_1}\norm{\hat{\A} - \A}), $$
which completes the proof of the theorem. $\square$

\vspace{12pt}

{\bf Proof of Theorem \ref{thm:statinvcov}.} We omit the rate for $\tilde{\bSigma}_{\y}$ since it involves standard treatments like that in Theorem \ref{thm:Invcov}. Also the proof for $\norm{\check{\bSigma}_{\y}^{-1} - \bSigma_{\y}^{-1}}$ follows exactly the arguments below for $\norm{\hat{\bSigma}_{\y}^{-1} - \bSigma_{\y}^{-1}}$, and is thus omitted.

Note that (\ref{eqn:thm:invcov:sigmayhatinv1}) becomes
$$ \norm{\hat{\bSigma}_{\y}^{-1} - \bSigma_{\y}^{-1}} = O_P((1+(p^{1-\delta_1}s^{-1})^{1/2})\norm{\hat{\A} - \A}) + O_P(L_1 + L_2), $$
with $L_1$ and $L_2$ defined similar to those in the proof of Theorem \ref{thm:Invcov}. Similar to (\ref{eqn:thm:invcov:L1}) and (\ref{eqn:thm:invcov:L2}), we have respectively
$$ L_1 = O_P(p^{-(1-\delta_2)}\norm{\hat{\A} - \A}), \;\;\; L_2 = O_P((1+(p^{1-\delta_1}s^{-1})^{1/2})\norm{\hat{\A} - \A}), $$
which shows that $ \norm{\hat{\bSigma}_{\y}^{-1} - \bSigma_{\y}^{-1}} = O_P((1+(p^{1-\delta_1}s^{-1})^{1/2})\norm{\hat{\A} - \A}) $. This completes the proof of the theorem. $\square$

\bibliographystyle{chicago}
\bibliography{myreference}
\end{document}